# SYMMETRIC SIMPLE EXCLUSION PROCESS IN DYNAMIC ENVIRONMENT: HYDRODYNAMICS

Frank Redig[1], Ellen Saada[2], Federico Sau[3]

April 28, 2020

### Abstract

We consider the symmetric simple exclusion process in $\mathbb{Z}^d$ with quenched bounded dynamic random conductances and prove its hydrodynamic limit in path space. The main tool is the connection, due to the self-duality of the process, between the invariance principle for single particles starting from all points and the macroscopic behavior of the density field. While the hydrodynamic limit at fixed macroscopic times is obtained via a generalization to the time-inhomogeneous context of the strategy introduced in [41], in order to prove tightness for the sequence of empirical density fields we develop a new criterion based on the notion of uniform conditional stochastic continuity, following [50]. In conclusion, we show that uniform elliptic dynamic conductances provide an example of environments in which the so-called arbitrary starting point invariance principle may be derived from the invariance principle of a single particle starting from the origin. Therefore, our hydrodynamics result applies to the examples of quenched environments considered in, e.g., [1], [3], [6] in combination with the hypothesis of uniform ellipticity.

## 1 Introduction

Dynamic random environments are natural quantities to be inserted in probabilistic models in order to make them more realistic. But studying such models is challenging, and for a long time only models endowed with a static environment were considered. However, *random walks in dynamic random environment* (RWDRE) have been extensively studied in recent years (see e.g. [1], [3], [4], [6], [7], [13], [45] and references therein) and several results on the law of large numbers, invariance principles and heat kernel estimates have been obtained. A natural next step is to consider particle systems in such dynamic environments. There the first question concerns the derivation of hydrodynamic limits. In this article, we answer this question for the nearest-neighbor symmetric simple exclusion process.

[1]Delft Institute of Applied Mathematics, TU Delft, Delft, The Netherlands. f.h.j.redig@tudelft.nl
[2]CNRS, UMR 8145, MAP5, Université de Paris, Campus Saint-Germain-des-Prés, Paris, France. ellen.saada@mi.parisdescartes.fr.
[3]Department of Mathematics & Computer Science, Institute of Science and Technology (IST Austria), 3400 Klosterneuburg, Austria. federico.sau@ist.ac.at

For interacting particle systems with a form of *self-duality* and that evolve in a *static disorder*, the problem of deriving the macroscopic equation governing the hydrodynamic limit has been shown to be strongly connected to the asymptotic behavior of a single random walker in the same environment. Indeed, the feature that if a rescaled test particle converges to a Brownian motion then the interacting particle system has a hydrodynamic limit appears already in e.g. [11], [24], [36] and [44]. Our contribution is to carry out this connection between single particle behavior and diffusive hydrodynamic limit in the context of *dynamic environment* for a nearest-neighbor particle system, namely the *symmetric simple exclusion process* (SSEP) in a *quenched dynamic bond disorder*, for which we show that a suitable form of self-duality remains valid. Let us now first recall the definition of SSEP, then detail the known results on SSEP evolving in a static environment.

**Symmetric simple exclusion process.** In words, the *symmetric simple exclusion process without disorder* in $\mathbb{Z}^d$ with $d \geq 1$ ([38], [47]) is an interacting particle system consisting of indistinguishable particles which are forbidden to simultaneously occupy the same site, and which jump at a constant rate only to nearest-neighbor unoccupied sites. More precisely, let $\eta \in \{0,1\}^{\mathbb{Z}^d}$ be a configuration of particles, with $\eta(x)$ denoting the number of particles at site $x \in \mathbb{Z}^d$. The stochastic process $\{\eta_t,\ t \geq 0\}$ is Markovian and evolves on the state space $\{0,1\}^{\mathbb{Z}^d}$ according to the infinitesimal generator

$$L\varphi(\eta) = \sum_{|x-y|=1} \left\{ \begin{array}{l} \eta(x)(1-\eta(y))(\varphi(\eta^{x,y}) - \varphi(\eta)) \\ + \eta(y)(1-\eta(x))(\varphi(\eta^{y,x}) - \varphi(\eta)) \end{array} \right\}, \tag{1.1}$$

where $|x-y| = \sum_{i=1}^{d} |x_i - y_i|$ and $\varphi : \{0,1\}^{\mathbb{Z}^d} \to \mathbb{R}$ is a bounded cylinder function, i.e. it depends only on a finite number of occupation variables $\{\eta(x),\ x \in \mathbb{Z}^d\}$. In (1.1) the finite summation is taken over all unordered pairs of nearest-neighboring sites – referred to as *bonds* – and $\eta^{x,y}$ is the configuration obtained from $\eta$ by removing a particle from the occupied site $x$ and placing it at the empty site $y$. The hydrodynamic limit ([11], [25], [34]) of the particle system described by (1.1) is known ([11], [34]) and, roughly speaking, prescribes that the trajectories of the particle density scale to the weak solution of the heat equation.

**Static environment.** For SSEP in a *quenched static bond disorder* in $\mathbb{Z}^d$, hydrodynamic limits – at a fixed macroscopic time – have been obtained by means of the *self-duality property* of the particle system, that is, by solving a homogenization problem (see e.g. [19, Theorem 2.1], [18, Theorem 2.4] and, more generally, [42]) or, alternatively, establishing an invariance principle (see e.g. [17], [41]) linked to the behavior of a single particle in the same environment. As examples, see [17], [41] for $d = 1$, [19] for $d \geq 1$ and [18] on the supercritical percolation cluster with $d \geq 2$. This method has been applied also to non-diffusive space-time rescalings, for which the hydrodynamic behavior is not described by a heat equation, see e.g. [17], [20], [21]. Nonetheless, all hydrodynamic limits obtained via this self-duality technique hold only at the level of finite-dimensional distributions and lack of a proof of relative compactness of the empirical density fields. Indeed, a direct application of the classical Aldous-Rebolledo criterion (see e.g. [34]) fails when following this approach.

Other techniques than self-duality – which apply to different particle systems – have also been studied in static environments. For instance, in quenched static bond disorder, the method based on the so-called *corrected empirical process* has been applied to prove hydrodynamics for SSEP ([31]) and for zero-range processes ([20], [27]). The non-gradient method ([44], [49], see also [34]) has found many applications to reversible lattice-gas models in a more general static environment, see e.g. [22].

**Dynamic environment.** In presence of *dynamic environment*, to the best of our knowledge, no hydrodynamic limit for interacting particle systems has been studied, yet.

On the one side, when looking at the hydrodynamic rescaling of a particle system in a quenched dynamic disorder, a space-time homogenization problem or, alternatively, an invariance principle for the associated RWDRE must be solved. On the other side, how this homogenization problem connects to the hydrodynamic behavior of the particle system depends on the interaction rules of the particles.

For the *symmetric simple exclusion process* in a *quenched dynamic bond disorder* in $\mathbb{Z}^d$ (whose generator is described in (2.3) below), we show that a form of *self-duality* still holds and allows us to write the occupation variables of the particle system in terms of positions of suitable *time-inhomogeneous backward random walks* evolving in the same environment. The hydrodynamic limit is, thus, obtained by studying the diffusive behavior of forward random walks evolving in this environment.

**A new tightness criterion.** In absence of criteria for relative compactness of the empirical density fields that apply to our case, we develop a tightness criterion based on the notion of *uniform conditional stochastic continuity* introduced in [50]. We formulate this tightness criterion to hold for a general sequence of $\mathbb{R}$-valued stochastic processes, though its validity extends straightforwardly to processes taking values in a general metric space. The main advantage of this criterion is that it avoids the use of stopping times as in Aldous criterion (see e.g. [32, Section 2.2]) to control the modulus of continuity of the processes' trajectories. In fact, a condition employing stopping times fails when the increments of the processes are expressed as stochastic integrals of non-predictable integrands, as in the case of stochastic convolutions as those in (3.7) below. However, to replace the "strong uniform stochastic equicontinuity" contained in Aldous criterion with a mere "uniform stochastic equicontinuity" as in [32, Remark 2.2.4] does not suffice to ensure tightness. With the purpose of bridging the gap between these two notions of stochastic equicontinuity, we show that a uniform control on the *conditional* tail probabilities of the processes' increments as in Theorem B.2 below suffices to guarantee tightness.

In this paper, we exploit this tightness criterion in two occasions. First, we use it in Appendix C to prove tightness of a sequence of random walks in a uniformly elliptic dynamic environment without making use of any estimates on hitting times of balls as done, e.g., in [10]. Secondly, we use it in Section 5.2 – combined with Mitoma's tightness criterion for tempered distribution-valued processes ([39]) – to prove tightness of the sequence of empirical density fields.

In the latter instance, we rely on two main assumptions for this tightness criterion to be effective: a quenched invariance principle for forward random walks and a uniform bound on the maximal number of particles per site.

In fact, under the aforementioned invariance principle hypothesis, this tightness criterion could be applied to other systems than the one considered in this paper. At first, because the (quenched) *static* bond disorder may be considered as a special instance of the *dynamic* environment we consider in this paper, the implementation of the tightness criterion to SSEP (see Proposition 5.5 below) carries directly through also for the particle systems in static environment considered e.g. in [17], [18], [19], [21], [41]. Furthermore, tightness may be proved via the same strategy for generalizations – in absence or presence of quenched static and/or dynamic bond disorder – of SSEP as in e.g. [9], in which up to $\alpha \in \mathbb{N}$ particles are allowed per site (note that this particle system differs from what is known as *generalized exclusion process*, see e.g. [34, Definition 2.4.1]). Even other types of quenched disorder are suited for this tightness criterion as proved in [23], in which the environment is designated by assigning a (uniformly bounded) maximal occupancy $\alpha_x \in \mathbb{N}$ to each site. In other words, this criterion applies to all particle systems for which a self-duality property and a uniform bound on the maximal number of particles per site hold in combination with the validity of the arbitrary starting point invariance principle.

**Arbitrary starting point invariance principle.** As we have already mentioned in the paragraph above, one of the two main assumptions is the validity of the invariance principle for random walks in a dynamic environment with arbitrary starting positions, equivalent, in turn, to semigroups and generators convergence (see, for more details, Theorem 5.2 below).

It is worth noticing that in the last ten years there have been several results in this direction, see e.g. [1], [3], [6], [7], [13], though all of them prove a quenched invariance principle for the "initially-anchored-at-the-origin" random walk only, i.e., for a given environment, the diffusively rescaled random walk that starts at time zero in the origin converges in law to a non-degenerate Brownian motion also starting in the origin. Unlike in the case of spatially homogeneous conductances, in our case the laws of the random walks are not translation invariant, therefore the derivation of an invariance principle for random walks centered around arbitrary macroscopic points does not follow at once from the invariance principle for the random walk initialized in the origin.

While this derivation cannot be proven to hold in general, some recent progress has been made in the case of static conductances in [10, Appendix A.2] and in the case of static site-inhomogeneities in [23, Section 4.4]. In Appendix C, we present the first instance involving dynamic environment in which this derivation from starting in the origin to arbitrary starting point invariance principle carries through. In particular, we establish this connection in the case of uniformly elliptic dynamic conductances providing, in combination with the results in e.g. [1], [3], [6], plenty of non-trivial examples of dynamic environments in which the assumption of arbitrary starting point invariance principle holds.

**Uniform boundedness of the environment.** Besides the validity of an arbitrary starting point invariance principle for the forward random walks, the other assumption on the environment that we require is its uniform boundedness – over the bonds and time. This assumption suffices to prove existence of the infinite particle system. Moreover, by means of this assumption alone and, in particular, without relying on uniform ellipticity of the environment (and consequent heat kernel estimates as in [1, Proposition 1.1]), we obtain an exponential upper bound for the transition probabilities of the random walks. This bound turns out to be useful in Appendix A.3 in providing an explicit formula for some observables of the particle system.

**Organization of the paper.** The rest of the paper is organized as follows. In Section 2 we introduce the dynamic environment and the model. In Section 3 we illustrate our approach in comparison with existing methods and state our main result, Theorem 3.2. In Section 4, from a graphical representation of the particle system (which we detail for the sake of completeness in Appendix A), we deduce a representation of the occupation variables as *mild solutions* of an infinite system of linear stochastic differential equations (which is proved in Appendix A.3). Section 5 is devoted to the proof of Theorem 3.2. We conclude the paper with the complete proof of our new tightness criterion used (Appendix B, more precisely Theorem B.4 in combination with Theorem B.2) and the study of a non-trivial space-time inhomogeneous scenario in which the invariance principle for the random walk starting from the origin yields an analogous invariance principle for all random walks starting from all macroscopic points and times (Appendix C, see also Section 3.1).

## 2 Setting

The space on which the particles move is the $d$-dimensional Euclidean lattice $(\mathbb{Z}^d, E_d)$ with $d \geq 1$. The set of bonds $E_d$ consists in all unordered pairs of nearest-neighboring sites, i.e.

$$E_d = \{\{x,y\}, \ x,y \in \mathbb{Z}^d \text{ with } |x-y| = 1\} \ .$$

Let us introduce our *dynamic environment* which is defined on the set of bonds $E_d$, so that we also refer to it as *(quenched) dynamic bond disorder* on $(\mathbb{Z}^d, E_d)$. Namely, we assign time-dependent non-negative weights to each bond $\{x,y\} \in E_d$ and we define as *environment* any càdlàg (w.r.t. the time variable $t$) function

$$\lambda = \{\lambda_t(\{x,y\}), \ \{x,y\} \in E_d, \ t \geq 0\} \ , \tag{2.1}$$

where
$$\lambda_t(\{x,y\}) = \lambda_t(\{y,x\}) \geq 0 \tag{2.2}$$
is referred to as the *conductance* of the bond $\{x,y\} \in E_d$ at time $t \geq 0$. The environment $\lambda$ is said to be *static* if $\lambda_t(\{x,y\}) = \lambda_0(\{x,y\})$ for all $\{x,y\} \in E_d$ and $t \geq 0$.

We will need the following assumption on the environment.

**Assumption 2.1** (BOUNDED CONDUCTANCES). *There exists a constant $\mathfrak{a} > 0$ for which, for all bonds $\{x,y\} \in E_d$ and $t \geq 0$, we have*
$$\lambda_t(\{x,y\}) \in [0, \mathfrak{a}].$$

**Remark 2.1.** *The boundedness of conductances guarantees, via a graphical construction (see Appendix A), that all stochastic processes introduced in Sections 3 and 4 are well-defined.*

Given the environment $\lambda$ as defined in (2.1)–(2.2), we now introduce as a counterpart to the symmetric simple exclusion process without disorder (1.1) the time-evolution of the *symmetric simple exclusion process in the dynamic environment $\lambda$* (SSEP($\lambda$)) by specifying its time-dependent infinitesimal generator $L_t$. For all $t \geq 0$ and every bounded cylinder function $\varphi : \{0,1\}^{\mathbb{Z}^d} \to \mathbb{R}$, we have

$$L_t \varphi(\eta) = \sum_{\{x,y\} \in E_d} \lambda_t(\{x,y\}) \left\{ \begin{array}{l} \eta(x)(1-\eta(y))(\varphi(\eta^{x,y}) - \varphi(\eta)) \\ + \eta(y)(1-\eta(x))(\varphi(\eta^{y,x}) - \varphi(\eta)) \end{array} \right\}. \tag{2.3}$$

Given any initial configuration $\eta \in \{0,1\}^{\mathbb{Z}^d}$, the time-dependent infinitesimal generators in (2.3) generate a time-inhomogeneous Markov (Feller) process $\{\eta_t, t \geq 0\}$ with sample paths in the Skorokhod space $D([0,\infty), \{0,1\}^{\mathbb{Z}^d})$ such that $\eta_0 = \eta$. We postpone to Section 4 the construction of this infinite particle system via a graphical representation.

## 3 Hydrodynamics result

In the present section we discuss the hydrodynamic limit in path space of the particle system $\{\eta_t, t \geq 0\}$ evolving in the environment $\lambda$, described by (2.3), that is, roughly speaking, the convergence in law of empirical density fields' trajectories to (deterministic) measures whose density w.r.t. the Lebesgue measure is solution of a Cauchy problem. Let us first detail what these density fields and the Cauchy problem with its solution are in our case.

**Empirical density fields.** We introduce for all $N \in \mathbb{N}$ the *empirical density field* $\{X_t^N, t \geq 0\}$ as a process in $D([0,\infty), \mathscr{S}'(\mathbb{R}^d))$, the Skorokhod space of $\mathscr{S}'(\mathbb{R}^d)$-valued càdlàg trajectories (see e.g. [39]), where $\mathscr{S}(\mathbb{R}^d)$ denotes the Schwartz class of rapidly decreasing functions on $\mathbb{R}^d$ and $\mathscr{S}'(\mathbb{R}^d)$ its topological dual. Given the particle system $\{\eta_t, t \geq 0\}$ evolving in the environment $\lambda$, for any test function $G \in \mathscr{S}(\mathbb{R}^d)$, the empirical density evaluated at $G$ reads as

$$\mathsf{X}_t^N(G) := \frac{1}{N^d} \sum_{x \in \mathbb{Z}^d} G(\tfrac{x}{N}) \eta_{tN^2}(x), \qquad t \geq 0. \tag{3.1}$$

So we choose to view the empirical density field as taking values in the space of tempered distributions rather than in the space of Radon measures as e.g. in [17], [19]. Indeed, the space $\mathscr{S}'(\mathbb{R}^d)$ has the advantage that it is a good space for tightness criteria (see e.g. [39]) and we use the fact that $\mathscr{S}(\mathbb{R}^d)$ is closed under the action of the Brownian motion semigroup.

**Heat equation.** Let $\langle \cdot, \cdot \rangle$ denote the standard scalar product in $\mathbb{R}^d$. We denote by $\{\rho_t^\Sigma, t \geq 0\}$ the *unique* weak solution to the following Cauchy problem

$$\begin{cases} \partial_t \rho &= \tfrac{1}{2} \nabla \cdot (\Sigma \nabla \rho) \\ \rho_0 &= \rho_\bullet, \end{cases} \tag{3.2}$$

with $\rho_\bullet : \mathbb{R}^d \to [0,1]$ measurable and $\Sigma \in \mathbb{R}^{d \times d}$ being a $d$-dimensional real symmetric positive-definite matrix (see e.g. [16], [34]). We recall that, for $\{\rho_t^\Sigma, t \geq 0\}$, being a weak solution of (3.2) means that, for all $G \in \mathscr{S}(\mathbb{R}^d)$ and $t \geq 0$,

$$\langle \rho_t^\Sigma, G \rangle = \langle \rho_\bullet, G \rangle + \int_0^t \langle \rho_s^\Sigma, \tfrac{1}{2} \nabla \cdot (\Sigma \nabla G) \rangle \, ds. \tag{3.3}$$

In addition, due to the linearity of the heat equation in (3.2), $\{\rho_t^\Sigma, t \geq 0\}$ may be represented in terms of $\{\mathcal{S}_t^\Sigma, t \geq 0\}$, the transition semigroup associated to the $d$-dimensional Brownian motion $\{B_t^\Sigma, t \geq 0\}$, starting at the origin and with covariance matrix $\Sigma \in \mathbb{R}^{d \times d}$ ([16]); namely, for all $G \in \mathscr{S}(\mathbb{R}^d)$,

$$\langle \rho_t^\Sigma, G \rangle = \langle \rho_\bullet, \mathcal{S}_t^\Sigma G \rangle. \tag{3.4}$$

**Hydrodynamics.** The proof of hydrodynamic limits in path space may be divided into two steps. First, one proves that, for all $T > 0$, the sequence of distributions of the empirical density fields $\{\mathsf{X}_t^N, 0 \leq t \leq T\}$ is relatively compact in $D([0,T], \mathscr{S}'(\mathbb{R}^d))$ by proving *tightness*. Then, one proves that all limiting probability measures in $D([0,T], \mathscr{S}'(\mathbb{R}^d))$ are supported on weak solutions of a Cauchy problem. By uniqueness of such a solution, the proof is concluded.

The "standard way" (e.g. [34]) to proceed is the following. To derive the convergence of the processes $\{\mathsf{X}_t^N, t \geq 0\}$ in $D([0,\infty), \mathscr{S}'(\mathbb{R}^d))$, we start from Dynkin's formula for the empirical density fields, i.e. for all $N \in \mathbb{N}$, $G \in \mathscr{S}(\mathbb{R}^d)$ and $t \geq 0$,

$$\mathsf{X}_t^N(G) = \mathsf{X}_0^N(G) + \int_0^t L_s \mathsf{X}_s^N(G) \, ds + \mathsf{M}_t^N(G), \tag{3.5}$$

with $\{\mathsf{M}_t^N(G), t \geq 0\}$ being a martingale. After obtaining tightness of the sequence

$\{X_t^N, \ 0 \leq t \leq T\}$ via an application of Aldous-Rebolledo criterion, the rest of the proof is carried out in two steps. First one shows that the martingale term $M_t^N(G)$ vanishes in probability as $N \to \infty$. Secondly, all the remaining terms in (3.5) can be expressed in terms of the empirical density field only; i.e. one "closes" the equation, yielding as a unique limit the solution expressed as in (3.3).

**Hydrodynamics & self-duality.** In presence of (static or dynamic) disorder, the issue of "closing" equation (3.5) in terms of the empirical density field only cannot be directly achieved. To overpass this obstacle, in the static disorder case, the authors in [27], [31] solve this problem by introducing an auxiliary observable, called *corrected empirical density field*.

Here we follow the probabilistic approach initiated in [41] and further developed in e.g. [17], which is more natural in our context. Key ingredients of this method are the *self-duality property* of the particle system and an alternative to Dynkin's formula (3.5), namely representing the empirical density fields as *mild solutions* of an infinite system of nonlinear – though with linear drift – stochastic differential equations: for all $N \in \mathbb{N}$, $G \in \mathscr{S}(\mathbb{R}^d)$ and $t \geq 0$,

$$X_t^N(G) = X_0^N(\mathcal{S}_{0,t}^N G) + \int_0^t dM_s^N(\mathcal{S}_{s,t}^N G), \qquad (3.6)$$

where, for all $s \in [0, \infty)$, $\{\mathcal{S}_{s,t}^N, \ t \in [s, \infty)\}$ may be related to the semigroup of a *suitably rescaled random walk* (see also (5.1) below) and with

$$\int_0^t dM_s^N(\mathcal{S}_{s,t}^N G) \qquad (3.7)$$

being typically indicated as a *stochastic convolution* term (see e.g. [43]).

Via this approach, the hydrodynamic limit is obtained in two steps: after proving tightness of the sequence of empirical density fields, first one ensures that the second term on the r.h.s. of (3.6) – which is not a martingale – vanishes in probability as $N \to \infty$; then, one checks whether the first term on the r.h.s. in (3.6) converges to $\langle \rho_t^\Sigma, G \rangle$ as given in (3.4), that is, $\langle \rho_t^\Sigma, G \rangle = \langle \rho_\bullet, \mathcal{S}_t^\Sigma G \rangle$. This latter convergence requires two ingredients: first, that the initial particle empirical density fields rescale (in probability) to a macroscopic density (assumption (a) in Theorem 3.2 below); secondly, that all random walks with arbitrary starting positions and evolving in the same dynamic environment rescale to Brownian motions with a given – space and time-independent – covariance matrix; namely, assumption (b) in Theorem 3.2 below.

In conclusion, while (3.3) is the representation of the solution $\{\rho_t^\Sigma, \ t \geq 0\}$ to the Cauchy problem (3.2) most commonly used when deriving hydrodynamic limits starting from Dynkin's formula (3.5), a method as the one we follow, based on the duality property of the particle system with suitable random walks, profits from the "mild solution" representation (3.4) of $\{\rho_t^\Sigma, \ t \geq 0\}$.

Let us now introduce the random walks alluded to above, used in our hydrodynamics result.

**Definition 3.1** (FORWARD AND BACKWARD RANDOM WALKS). *For all $s \geq 0$, let $\{X^x_{s,t},\ t \in [s, \infty)\}$ be the* forward *random walk starting at $x \in \mathbb{Z}^d$ at time $s$ and evolving in the environment $\lambda$ through the time-dependent infinitesimal generator*

$$A_t f(x) = \sum_{y:\{x,y\} \in E_d} \lambda_t(\{x,y\})(f(y) - f(x)),\tag{3.8}$$

*where $f : \mathbb{Z}^d \to \mathbb{R}$ is a bounded function.*

*Similarly, for all $t \geq 0$, let $\{\widehat{X}^y_{s,t},\ s \in [0,t]\}$ be the* backward *random walk which starts at $y \in \mathbb{Z}^d$ at time $t$ and "evolves backwards" in the environment $\lambda$ through the time-dependent infinitesimal generator*

$$A_{s^-} f(x) = \sum_{y:\{x,y\} \in E_d} \lambda_{s^-}(\{x,y\})(f(y) - f(x)),\tag{3.9}$$

*where $f : \mathbb{Z}^d \to \mathbb{R}$ is as above and $\lambda_{s^-}(\{x,y\}) = \lim_{r \uparrow s} \lambda_r(\{x,y\})$ for all $s \in [0,t]$.*

We will give in Section 4.1 and Appendix A.1 the construction of both those forward and backward random walks via a graphical representation.

We are now ready to state our main theorem, Theorem 3.2, followed by two remarks related to its proof.

**Theorem 3.2** (PATH-SPACE HYDRODYNAMIC LIMIT). *For all $N \in \mathbb{N}$, we initialize the exclusion process $\{\eta_t,\ t \geq 0\}$ according to a probability measure $\mu_N$ on $\{0,1\}^{\mathbb{Z}^d}$ (Notation: $\eta_0 \sim \mu_N$), and, consequently, $\mathsf{X}^N_0$ is the random element in $\mathscr{S}'(\mathbb{R}^d)$ obtained from $\eta_0 \sim \mu_N$. Besides Assumption 2.1, we further assume that*

(a) *The family of probability measures $\{\mu_N,\ N \in \mathbb{N}\}$ on $\{0,1\}^{\mathbb{Z}^d}$ is associated to the density profile $\rho_\bullet : \mathbb{R}^d \to [0,1]$; namely, for all $G \in \mathscr{S}(\mathbb{R}^d)$ and $\delta > 0$,*

$$\mu_N\left(\left|\frac{1}{N^d}\sum_{x \in \mathbb{Z}^d} G(\tfrac{x}{N})\eta(x) - \int_{\mathbb{R}^d} G(u)\rho_\bullet(u)\,\mathrm{d}u\right| > \delta\right) \xrightarrow[N \to \infty]{} 0.\tag{3.10}$$

(b) *The forward random walks $\{X^x_{0,t},\ x \in \mathbb{Z}^d,\ t \in [0,\infty)\}$ with arbitrary starting positions satisfy an* invariance principle *with a non-degenerate covariance matrix $\Sigma$; namely, for all $u \in \mathbb{R}^d$, for all sequences $\{x_N,\ N \in \mathbb{N}\} \subset \mathbb{Z}^d$ for which $\frac{x_N}{N} \to u \in \mathbb{R}^d$ as $N \to \infty$ and for all $T > 0$,*

$$\left\{\frac{X^{x_N}_{0,tN^2}}{N},\ t \in [0,T]\right\} \underset{N \to \infty}{\Longrightarrow} \{B^\Sigma_t + u,\ t \in [0,T]\},\tag{3.11}$$

*where $\{B_t^\Sigma, \ t \geq 0\}$ denotes the d-dimensional Brownian motion introduced below (3.3), starting at the origin and with non-degenerate covariance matrix $\Sigma$ (Notation: $\Rightarrow$ stands for convergence in law).*

*Then, for all $T > 0$, we have the following convergence*

$$\{X_t^N, \ t \in [0, T]\} \underset{N \to \infty}{\Longrightarrow} \{\pi_t^\Sigma, \ t \in [0, T]\} \tag{3.12}$$

*in $D([0, T], \mathscr{S}'(\mathbb{R}^d))$, with $\pi_t^\Sigma(du) = \rho_t^\Sigma(u)\,du$ and $\{\rho_t^\Sigma, \ t \geq 0\}$ being the unique weak solution to the Cauchy problem (3.2).*

**Remark 3.3** (UNIFORM CONVERGENCE OVER TIME). *If $\pi_t^\Sigma(du) = \rho_t^\Sigma(u)du$ for all $t \geq 0$ and $\rho_0 = \rho_\bullet \in L^\infty(\mathbb{R}^d)$ (which holds true in our setting), then $\{\pi_t^\Sigma, \ t \geq 0\}$ is a trajectory in $C([0, \infty), \mathscr{S}'(\mathbb{R}^d))$, the space of tempered distribution-valued continuous trajectories. Indeed, for all $G \in \mathscr{S}(\mathbb{R}^d)$ and $t \geq 0$, as long as $\rho_\bullet \in L^\infty(\mathbb{R}^d)$, by (3.4), we have*

$$\left|\langle \rho_t^\Sigma, G\rangle - \langle \rho_s^\Sigma, G\rangle\right| \leq \left(\sup_{u \in \mathbb{R}^d} \rho_\bullet(u)\right) \int_{\mathbb{R}^d} \left|\mathcal{S}_t^\Sigma G(u) - \mathcal{S}_s^\Sigma G(u)\right| du \underset{|t-s| \to 0}{\longrightarrow} 0\,. \tag{3.13}$$

*Hence, because weakly continuous trajectories in the space $\mathscr{S}'(\mathbb{R}^d)$ are strongly continuous (see e.g. [30, p. 145]), it follows that $\{\pi_t^\Sigma, \ t \geq 0\} \in C([0, \infty), \mathscr{S}'(\mathbb{R}^d))$.*

*As a consequence, the convergence in $D([0, T], \mathscr{S}'(\mathbb{R}^d))$ in (3.12) becomes convergence w.r.t. the uniform topology of $C([0, T], \mathscr{S}'(\mathbb{R}^d))$, i.e. it can also be equivalently rewritten as follows: for all $G \in \mathscr{S}(\mathbb{R}^d)$,*

$$\sup_{0 \leq t \leq T} \left|\frac{1}{N^d} \sum_{x \in \mathbb{Z}^d} G(\tfrac{x}{N})\eta_{tN^2}(x) - \int_{\mathbb{R}^d} G(u)\, \rho_t^\Sigma(u)\, du\right| \underset{N \to \infty}{\overset{\mathcal{P}}{\longrightarrow}} 0\,,$$

*(Notation: $\overset{\mathcal{P}}{\to}$ stands for convergence in probability). Indeed, by e.g. [39, Theorem 5.3.2], for all $G \in \mathscr{S}(\mathbb{R}^d)$, the real-valued stochastic processes $\{X_t^N(G), \ t \in [0, T]\}$ in $D([0, T], \mathbb{R})$ converge in law (w.r.t. the Skorokhod topology of $D([0, T], \mathbb{R})$) to the (deterministic) process $\{\langle \rho_t^\Sigma, G\rangle, \ t \in [0, T]\}$. On the other hand, because $\{\langle \rho_t^\Sigma, G\rangle, \ t \in [0, T]\}$ belongs to $C([0, T], \mathbb{R})$, by e.g. [5, end of p. 124], the processes $\{X_t^N(G), \ t \in [0, T]\}$ converge in law – and, thus, in probability – to the same deterministic limit w.r.t. the uniform topology of $C([0, T], \mathbb{R})$.*

We will prove Theorem 3.2 in Section 5. As explained earlier, it relies on a mild solution representation for the particle system involving the forward and backward random walks of Definition 3.1 (see Section 4). This representation induces a mild solution representation of the empirical density fields, obtained in formula (5.1).

The proof of tightness for the empirical density fields – which cannot be achieved by means of more standard techniques (e.g. Censov or Aldous-Rebolledo criteria, resp. to be found e.g. in [11] and [34]) when representing the fields in this form – requires the elaboration of a new tightness criterion – Theorem B.4 in combination with Theorem B.2

–, presented in Appendix B. The proof of relative compactness, carried out in Section 5.2, will use this criterion, formula (5.1), and the arbitrary starting point invariance principle for the forward random walks.

The characterization of the limiting measures as concentrated on the unique weak solution of the Cauchy problem (3.2) boils down to prove convergence of finite-dimensional distributions, that is: for all $n \in \mathbb{N}$, for all $0 \leq t_1 < \ldots < t_n \leq T$ and for all $G_1, \ldots, G_n \in \mathscr{S}(\mathbb{R}^d)$,

$$\left(X_{t_1}^N(G_1), \ldots, X_{t_n}^N(G_n)\right) \xrightarrow[N \to \infty]{\mathcal{P}} \left(\langle \rho_{t_1}^{\Sigma}, G_1\rangle, \ldots, \langle \rho_{t_n}^{\Sigma}, G_n\rangle\right) . \tag{3.14}$$

As joint convergence in probability comes down to checking convergence in probability of the single marginal laws, it suffices to prove (3.14) for the choice $n = 1$ only, that is: for all $0 \leq t \leq T$ and $G \in \mathscr{S}(\mathbb{R}^d)$,

$$\left| \frac{1}{N^d} \sum_{x \in \mathbb{Z}^d} G(\tfrac{x}{N}) \eta_{tN^2}(x) - \int_{\mathbb{R}^d} G(u) \rho_t^{\Sigma}(u) \, du \right| \xrightarrow[N \to \infty]{\mathcal{P}} 0 . \tag{3.15}$$

In Section 5.1, we will then exploit the mild solution representation (5.1) to prove (3.15). For this, we will further generalize to the time-inhomogeneous context results originally developed in [41] and further extended in e.g. [17], [18], [19].

We end this section with a discussion on examples of dynamic environments yielding condition (b) in Theorem 3.2.

## 3.1 Condition (b) & examples of dynamic environments

Our assumption (b) in Theorem 3.2 may be seen as the dynamic counterpart of the static arbitrary starting point quenched invariance principle in [41, Theorem 1] and [17, Proposition 4.3]. There, both authors derive this crucial result – rather than assuming it, as we do – from statistical properties of the conductances, namely strictly positive and uniformly bounded i.i.d. conductances with a fourth-negative moment condition and a (strong) law of large numbers for the inverse of the conductances (= resistances), respectively. Via different techniques, the same authors show that those two assumptions suffice – in dimension $d = 1$ and in presence of static conductances – to derive [41, Theorem 1] and [17, Proposition 4.3], which, by Theorem 5.2 below, are both equivalent to our assumption (b) in presence of static environment.

In recent years (see e.g. [1], [3], [4], [6], [7], [13], [45]) there has been an intensive research in providing general examples of dynamical environments $\lambda$ leading to non-degenerate invariance principles for the forward random walk $\{X_{0,t}^0, t \geq 0\}$ starting at the origin $0 \in \mathbb{Z}^d$. In all these cases, $\lambda$ is obtained as a typical realization of a suitably constructed random environment process $(\Lambda, F, P)$, yielding, for $P$-a.e. environment $\lambda \in$

$\Lambda$ and all $T > 0$,

$$\left\{\frac{X^0_{0,tN^2}}{N}, \ t \in [0,T]\right\} \underset{N \to \infty}{\Longrightarrow} \left\{B^\Sigma_t, \ t \in [0,T]\right\}. \tag{3.16}$$

Several examples of dynamic random environments which lead to invariance principles as those in (3.16) have been studied in the aforementioned references. In particular, it is worth mentioning that dynamic random environments driven by i.i.d. flipping and Markovian conductances taking values on a finite subset of $(0, \infty)$ fall in the setting studied in [1] for all dimensions $d \geq 1$. Further examples are provided in [7], [40]. There, the authors consider – among other examples – the symmetric simple exclusion process in $\mathbb{Z}^d$ with $d \geq 2$ as an interacting particle system which induces the underlying dynamic random environment for the random walk.

In fact, more general random environments that fit our context have been studied. In particular, in the works [1], [3], [6], [7], [13], the authors obtain quenched invariance principles with deterministic and non-degenerate covariance matrices $\Sigma$ for space-time ergodic random dynamical environments under conditions of either ellipticity or boundedness of $p$-moments of conductances and resistances.

However, all these quenched invariance principles in dynamic random environment are obtained for the random walk initially anchored at the origin, whilst our assumption (b) in Theorem 3.2 consists in a quenched invariance principle holding simultaneously for *all* random walks centered around *all* macroscopic points $u \in \mathbb{R}^d$. The problem of deriving such "arbitrary starting point quenched invariance principles" (our assumption (b)) from (3.16) has been addressed in the case of static environment in [10, Appendix A.2], while – to the best of our knowledge – it has remained unsolved in the dynamic setting.

To this purpose, we derive in Theorem 3.4 below – under the stronger assumption of uniformly elliptic dynamic environment, i.e. assumption ($b_2$) in the same theorem – such arbitrary starting point invariance principle, i.e. condition (b) in Theorem 3.2, from the invariance principle for the random walk starting at the origin, i.e. assumption ($b_1$) in Theorem 3.4 below. We postpone the proof of this theorem to Appendix C.

**Theorem 3.4** (ARBITRARY STARTING POINT INVARIANCE PRINCIPLE). *We assume that:*

(b$_1$) *The invariance principle for the random walk starting in the origin $0 \in \mathbb{Z}^d$ at time $s = 0$ holds; namely, for all $T > 0$, (3.16) holds.*

(b$_2$) *There exist constants $\mathfrak{a} > \mathfrak{b} > 0$ for which, for all bonds $\{x,y\} \in E_d$ and $t \geq 0$, we have*

$$\lambda_t(\{x,y\}) \in [\mathfrak{b}, \mathfrak{a}].$$

*Then, the invariance principle with arbitrary starting points (condition* (b) *of Theorem 3.2 ) holds; namely, for all $T > 0$, for all $u \in \mathbb{R}^d$ (macroscopic points) and for all $\{x_N : N \in \mathbb{N}\} \subset \mathbb{Z}^d$ such that $\frac{x_N}{N} \to u$ as $N \to \infty$ (approximating points on the rescaled lattice), (3.11) holds.*

# 4 Graphical constructions and mild solution

In Section 4.1 we construct the symmetric simple exclusion process in dynamic environment via a graphical representation. Relying on this construction, we express in Section 4.2 the occupation variables of the symmetric simple exclusion process (viewed as a stirring process) in dynamic environment as mild solution of a system of Poissonian stochastic differential equations.

## 4.1 Graphical construction of the particle system

The graphical construction employs, as a source of randomness, a collection of independent Poisson processes, each one attached to a bond of $\mathbb{Z}^d$. To take care of both space and time inhomogeneities, their intensities will depend both on the bond and time. As an intermediate step towards the graphical construction of the particle system, the same Poisson processes provide a graphical construction for all forward and backward random walks introduced in Definition 3.1. We explain this procedure below, leaving a detailed treatment to Appendix A. Finally, we will relate the occupation variables of the particle system to the positions of backward random walks. This must be meant in a pathwise sense, expressing the pathwise duality of the symmetric simple exclusion process in the dynamic environment $\lambda$.

**Poisson processes.** We consider a family of independent inhomogeneous Poisson processes

$$\{\mathcal{N}.(\{x,y\}), \{x,y\} \in E_d\} \tag{4.1}$$

defined on the probability space $(\Omega, \mathcal{F}, \{\mathcal{F}_t, t \geq 0\}, \mathbb{P})$, where $\{\mathcal{F}_t, t \geq 0\}$ is the natural filtration, $\mathcal{F} = \sigma\left(\cup_{t \geq 0} \mathcal{F}_t\right)$ and such that $\mathcal{N}.(\{x,y\})$ has intensity measure $\lambda_r(\{x,y\}) \, \mathrm{d}r$, that is

$$\mathbb{E}[\mathcal{N}_t(\{x,y\})] = \int_0^t \lambda_r(\{x,y\}) \, \mathrm{d}r, \qquad t \geq 0,$$

where $\mathbb{E}$ denotes expectation w.r.t. $\mathbb{P}$ (for a constructive definition of the probability space $(\Omega, \mathcal{F}, \{\mathcal{F}_t, t \geq 0\}, \mathbb{P})$, we refer to Appendix A.1). The associated compensated Poisson processes

$$\{\overline{\mathcal{N}}.(\{x,y\}), \{x,y\} \in E_d\},$$

defined as

$$\overline{\mathcal{N}}_t(\{x,y\}) = \mathcal{N}_t(\{x,y\}) - \int_0^t \lambda_r(\{x,y\}) \, \mathrm{d}r, \quad t \geq 0, \tag{4.2}$$

are a family of square integrable martingales w.r.t. $\{\mathcal{F}_t, t \geq 0\}$ of bounded variation, due to Assumption 2.1.

The associated picture is drawn as follows. On the space $\mathbb{Z}^d \times [0, \infty)$, where $\mathbb{Z}^d$ represents the sites and $[0, \infty)$ represents time which goes up, for each $z \in \mathbb{Z}^d$ draw a vertical line $\{z\} \times [0, \infty)$. Then for each $\{x,y\} \in E_d$, draw a horizontal two-sided arrow between $x$ and $y$ at each event time, i.e. jump time, of $\mathcal{N}.(\{x,y\})$.

**Forward and backward random walks.** We recover the walks defined in Definition 3.1 as follows. First, for all $\omega \in \Omega$, $x \in \mathbb{Z}^d$, $s \geq 0$ and $t \geq s$, $X^x_{s,t}[\omega]$ now denotes the position at time $t$ of the random walk in $\mathbb{Z}^d$ that is at $x$ at time $s$ and that, between times $s$ and $t$, crosses the bond $\{z, v\} \in E_d$ at an event time of $\mathcal{N}.(\{z, v\})[\omega]$ whenever at that time the walk is at location either $z$ or $v$ in $\mathbb{Z}^d$ (i.e. it follows the corresponding arrow in the graphical representation). We prove in Appendix A, thanks to Assumption 2.1, that the trajectories of those walks are, for $\mathbb{P}$-a.e. realization $\omega \in \Omega$, well defined for all times and starting positions. In fact, they are all simultaneously defined on the probability space $(\Omega, \mathcal{F}, \{\mathcal{F}_t, t \geq 0\}, \mathbb{P})$. In Appendix A, we show that their associated generators are given by (3.8), so that, indeed, these walks are a version of the processes introduced in Definition 3.1.

We now provide a version of the backward random walks of Definition 3.1. For $\mathbb{P}$-a.e. $\omega \in \Omega$, $t \geq 0$ and $y \in \mathbb{Z}^d$, we implicitly define *backward* random walks' trajectories $\{\widehat{X}^y_{s,t}[\omega], s \in [0, t]\}$ by the following identity:

$$X^{\widehat{X}^y_{s,t}[\omega]}_{s,t}[\omega] = y. \tag{4.3}$$

In words, $\widehat{X}^y_{s,t}[\omega]$ denotes the position in $\mathbb{Z}^d$ at time $s$ of the forward random walk that follows the Poissonian marks associated to $\omega \in \Omega$ and that is at $y \in \mathbb{Z}^d$ at time $t$ with $t \geq s$. In particular, for $\mathbb{P}$-a.e. $\omega \in \Omega$ and $x, y \in \mathbb{Z}^d$, we have

$$X^x_{s,t}[\omega] = y \quad \text{if and only if} \quad \widehat{X}^y_{s,t}[\omega] = x. \tag{4.4}$$

Again, all these random walks are simultaneously $\mathbb{P}$-a.s. well-defined, and these backward random walks coincide in law with the ones in Definition 3.1 (see Appendix A).

**Transition probabilities.** The Poissonian construction and the jump rules explained above ensure that each of the forward and backward random walks is Markovian.

For all $x, y \in \mathbb{Z}^d$, $s \geq 0$ and $t \geq s$, if we define

$$p_{s,t}(x, y) = \mathbb{P}\left(X^x_{s,t} = y\right) \quad \text{and} \quad \widehat{p}_{s,t}(y, x) = \mathbb{P}\left(\widehat{X}^y_{s,t} = x\right), \tag{4.5}$$

we obtain families of transition probabilities respectively for the forward and backward random walks. In particular, for all $x, y \in \mathbb{Z}^d$ and $0 \leq s \leq r \leq t$, we have the Chapman-Kolmogorov equations

$$\sum_{z \in \mathbb{Z}^d} p_{s,r}(x, z) p_{r,t}(z, y) = p_{s,t}(x, y) \tag{4.6}$$

$$\sum_{z \in \mathbb{Z}^d} \widehat{p}_{r,t}(y, z) \widehat{p}_{s,r}(z, x) = \widehat{p}_{s,t}(y, x). \tag{4.7}$$

From (4.4), we obtain that

$$p_{s,t}(x, y) = \widehat{p}_{s,t}(y, x), \tag{4.8}$$

for all $x, y \in \mathbb{Z}^d$ and $t \geq s$. Then, the operators $\{S_{s,t},\ t \in [s, \infty)\}$ and $\{\widehat{S}_{s,t},\ s \in [0, t]\}$, acting on bounded functions $f : \mathbb{Z}^d \to \mathbb{R}$ as, for $x \in \mathbb{Z}^d$,

$$S_{s,t} f(x) = \sum_{y \in \mathbb{Z}^d} p_{s,t}(x, y) f(y) \tag{4.9}$$

$$\widehat{S}_{s,t} f(x) = \sum_{y \in \mathbb{Z}^d} \widehat{p}_{s,t}(x, y) f(y), \tag{4.10}$$

correspond to the transition semigroups (or, more properly, the "evolution systems" or "forward/backward propagators" as referred to in [8] and references therein) respectively associated to the forward and backward random walks. Then, as a consequence of (4.8), we obtain that

$$\sum_{x \in \mathbb{Z}^d} [S_{s,t} f(x)] g(x) = \sum_{x \in \mathbb{Z}^d} f(x) \widehat{S}_{s,t} g(x), \tag{4.11}$$

for all $f, g : \mathbb{Z}^d \to \mathbb{R}$ for which the above summations are finite.

We refer to Appendix A.2 for a more detailed treatment with further properties of the above transition probabilities and associated time-inhomogeneous semigroups.

**Stirring process.** The stirring process relates the above introduced random walks with the occupation variables of the symmetric simple exclusion process in the environment $\lambda$ as follows. Due to the symmetry (2.2) of the environment and the one of the exclusion dynamics, we can rewrite the generator (2.3) as

$$L_t \varphi(\eta) = \sum_{\{x,y\} \in E_d} \lambda_t(\{x, y\}) (\varphi(\eta^{\{x,y\}}) - \varphi(\eta)),$$

where $\eta^{\{x,y\}}$ stands for the exchange of occupation numbers between sites $x$ and $y$ in configuration $\eta$, which takes place even if $x, y$ are both occupied (due to the fact that particles are indistinguishable). This rewriting gives the stirring interpretation of the symmetric simple exclusion process in the environment $\lambda$ (similar to the stirring interpretation in the case (1.1) without disorder, as described in [11, p. 98] and [38, p. 399]), that we take from now on. This way, the stirring process can be constructed on the same graphical representation as before, and particles evolve as the forward random walks previously introduced. Indeed, on $\mathbb{Z}^d \times [0, \infty)$, place a particle at $\{z\} \times \{0\}$ whenever $\eta(z) = 1$. Then the particle at $x$, if there is one, goes up on $\mathbb{Z}^d \times [0, \infty)$, following the random walk $X_{0,\cdot}^x[\omega]$.

Hence, similarly to [38, p. 399], we can write, for $\mathbb{P}$-a.e. $\omega \in \Omega$, for any initial configuration $\eta \in \{0, 1\}^{\mathbb{Z}^d}$, for any $x \in \mathbb{Z}^d$ and $t \geq 0$, that

$\eta_t(x)[\omega] = 1$ if and only if there is a $y \in \mathbb{Z}^d$ so that $X_{0,t}^y[\omega] = x$ and $\eta(y) = 1$

or, equivalently by using the associated backward random walks and (4.3),

$\eta_t(x)[\omega] = 1$ if and only if there is a $y \in \mathbb{Z}^d$ so that $\widehat{X}_{0,t}^x[\omega] = y$ and $\eta(y) = 1$ .

In other words,

$$\eta_t(x)[\omega] = \eta(\widehat{X}_{0,t}^x[\omega]), \qquad x \in \mathbb{Z}^d, \quad t \geq 0, \tag{4.12}$$

thus the stochastic process $\{\eta_t, t \geq 0\}$ (with $\eta_0 = \eta$) is defined for $\mathbb{P}$-a.e. $\omega \in \Omega$ on the probability space $(\Omega, \mathcal{F}, \{\mathcal{F}_t, t \geq 0\}, \mathbb{P})$. Moreover, from the memoryless property of the inhomogeneous Poisson processes employed in the graphical construction of forward and backward random walks, given any initial configuration $\eta \in \{0,1\}^{\mathbb{Z}^d}$, we recover the Markov property of the process $\{\eta_t, t \geq 0\}$ w.r.t. $\{\mathcal{F}_t, t \geq 0\}$.

**Remark 4.1** (PATHWISE SELF-DUALITY OF SSEP IN DYNAMIC ENVIRONMENT). *What we obtained in (4.12) is the property of pathwise self-duality of the symmetric simple exclusion process with a single dual particle (=a one-particle system backward in time), which thus remains valid also in presence of the dynamic environment $\lambda$.*

**Remark 4.2** (NOTATION). *In Theorem 3.2, we have $\eta_0 \sim \mu_N$. We thus have to enlarge $\Omega$ and, accordingly, the filtration and the probability measure, to simultaneously take into account possibly different initial particle configurations. Nevertheless, for the sake of simplicity, we will always write $(\Omega, \mathcal{F}, \{\mathcal{F}_t, t \geq 0\})$, but we will write $\mathbb{P}_{\mu_N}$ (resp. $\mathbb{P}_\eta$) for the probability measure induced by the Poisson processes in (4.1) and the distribution $\mu_N$ (resp. the Dirac measure $\delta_\eta$) of the initial configuration $\eta_0 \in \{0,1\}^{\mathbb{Z}^d}$ of the exclusion process $\{\eta_t, t \geq 0\}$ (and $\mathbb{E}_{\mu_N}$ (resp. $\mathbb{E}_\eta$) for the corresponding expectation).*

## 4.2 Mild solution representation of the particle system

The above construction provides an alternative way of defining $\{\eta_t, t \geq 0\}$, the symmetric simple exclusion process in the environment $\lambda$ as strong solution of an infinite system of linear stochastic differential equations. This is the content of Proposition 4.3 below. For an analogous statement previously obtained in the time-homogeneous context, we refer to identity (13) in [41].

The motivation comes from an infinitesimal description of the stirring process, as explained through the following computation. For all $t > 0$ and $x \in \mathbb{Z}^d$, if we write $d\eta_t(x) = \eta_t(x) - \eta_{t^-}(x)$, we have

$$d\eta_t(x) = \sum_{y:\{x,y\}\in E_d} (\eta_{t^-}(y) - \eta_{t^-}(x)) \, d\mathcal{N}_t(\{x,y\}). \tag{4.13}$$

By introducing the compensated Poisson process (4.2) in (4.13), we obtain

$$d\eta_t(x) = \sum_{y:\{x,y\}\in E_d} (\eta_{t^-}(y) - \eta_{t^-}(x)) \lambda_t(\{x,y\}) \, dt \\ + \sum_{y:\{x,y\}\in E_d} (\eta_{t^-}(y) - \eta_{t^-}(x)) \, d\overline{\mathcal{N}}_t(\{x,y\}). \tag{4.14}$$

Note that the terms in the second sum in the r.h.s. of (4.14) are increments of a martingale as products of bounded predictable terms and increments of the compensated

Poisson processes. Moreover, like the latter, such martingales are square integrable and of bounded variation.

After observing that the first sum on the r.h.s. of (4.14) corresponds to the definition of the infinitesimal generator in (3.8) at time $t$ of the forward random walk, we rewrite (4.14) as

$$d\eta_t(x) = A_t\,\eta_{t^-}(x)\,dt + dM_t(\eta_{t^-}, x)\,, \qquad x \in \mathbb{Z}^d\,,\quad t > 0\,, \qquad (4.15)$$

where $A_t$ acts on the $x$ variable and where

$$dM_t(\eta, x) := \sum_{y:\{x,y\}\in E_d} (\eta(y) - \eta(x))\,d\overline{\mathcal{N}}_t(\{x,y\})\,. \qquad (4.16)$$

In the following proposition, whose proof is postponed to Appendix A.3, we state that the so-called "mild solution" [43, Chapter 9] associated to the system of differential equations (4.15) equals $\mathbb{P}$-a.s. the process obtained via the stirring procedure in (4.12). The mild solution is defined as in (4.17) below, i.e. by formally applying the method of variation of constants to (4.15). Recall that $\{\widehat{S}_{s,t},\ s \in [0,t]\}$ and $\{\widehat{p}_{s,t}(y,x),\ x,y \in \mathbb{Z}^d,\ s \in [0,t]\}$ are, respectively, the semigroup and transition probabilities (see (4.10) and (4.5), respectively) of the backward random walks of Definition 3.1.

**Proposition 4.3.** *Fix an initial configuration $\eta \in \{0,1\}^{\mathbb{Z}^d}$. Consider, for all $x \in \mathbb{Z}^d$ and $t \geq 0$,*

$$\zeta_t(x) = \widehat{S}_{0,t}\eta(x) + \int_0^t \widehat{S}_{r,t}\,dM_r(\eta_{r^-}, x)$$

$$= \sum_{y\in\mathbb{Z}^d} \widehat{p}_{0,t}(x,y)\eta(y) + \int_0^t \sum_{y\in\mathbb{Z}^d} \widehat{p}_{r,t}(x,y)\,dM_r(\eta_{r^-}, y)\,, \qquad (4.17)$$

*where $\{\eta_t(x),\ x \in \mathbb{Z}^d,\ t \geq 0\}$ is defined in (4.12), $\eta_0 = \eta$ and $dM_r$ is given in (4.16). Then, for $\mathbb{P}$-a.e. $\omega \in \Omega$,*

$$\zeta_t(x)[\omega] = \eta_t(x)[\omega]\,, \qquad x \in \mathbb{Z}^d\,,\quad t \geq 0\,.$$

**Remark 4.4.** *Systems of equations of type (4.15) are studied in [43] in the context of Hilbert spaces. There it is proved that for a large class of semi-linear infinite-dimensional SDEs the so-called mild solutions coincide with weak solutions.*

## 5 Proof of Theorem 3.2

The key ingredient to prove Theorem 3.2 is the decomposition of the occupation variables of the process $\{\eta_t,\ t \geq 0\}$ provided in Proposition 4.3.

Let $G \in \mathscr{S}(\mathbb{R}^d)$ and $\eta_0 = \eta \in \{0,1\}^{\mathbb{Z}^d}$ be a fixed initial configuration. We consider the empirical density fields $X_t^N(G)$ as defined in (3.1). By using Proposition 4.3 and, then, identity (4.11), we obtain, by viewing the configuration $\eta$ as a function $\eta : \mathbb{Z}^d \to \{0, 1\}$ of the random walk's position,

$$\begin{aligned}
\mathsf{X}_t^N(G) &= \frac{1}{N^d} \sum_{x \in \mathbb{Z}^d} G(\tfrac{x}{N}) \eta_{tN^2}(x) \\
&= \frac{1}{N^d} \sum_{x \in \mathbb{Z}^d} G(\tfrac{x}{N}) \widehat{S}_{0,tN^2} \eta(x) + \frac{1}{N^d} \sum_{x \in \mathbb{Z}^d} G(\tfrac{x}{N}) \int_0^{tN^2} \widehat{S}_{r,tN^2} \, \mathrm{d} M_r(\eta_{r^-}, x) \\
&= \frac{1}{N^d} \sum_{x \in \mathbb{Z}^d} [S_{0,tN^2}^N G(\tfrac{x}{N})] \eta(x) + \frac{1}{N^d} \sum_{x \in \mathbb{Z}^d} \int_0^{tN^2} [S_{r,tN^2}^N G(\tfrac{x}{N})] \, \mathrm{d} M_r(\eta_{r^-}, x) \\
&= \mathsf{X}_0^N(S_{0,tN^2}^N G) + \frac{1}{N^d} \sum_{x \in \mathbb{Z}^d} \int_0^{tN^2} [S_{r,tN^2}^N G(\tfrac{x}{N})] \, \mathrm{d} M_r(\eta_{r^-}, x), \quad (5.1)
\end{aligned}$$

where

$$S_{s,t}^N G(\tfrac{x}{N}) := S_{s,t} G(\tfrac{\cdot}{N})(x), \qquad x \in \mathbb{Z}^d. \quad (5.2)$$

Note that the decomposition (5.1) (different from Dynkin's formula (3.5)) is the one presented in (3.6).

We then proceed as announced after (3.6): in Section 5.2 we exploit the tightness criterion given in Appendix B to prove relative compactness of the empirical density fields. In Section 5.1 we prove convergence of finite-dimensional distributions, that is (3.15), by showing that, for any $\delta > 0$,

$$\mathbb{P}_{\mu_N}\left(\left| \frac{1}{N^d} \sum_{x \in \mathbb{Z}^d} \int_0^{tN^2} S_{r,tN^2}^N G(\tfrac{x}{N}) \, \mathrm{d} M_r(\eta_{r^-}, x) \right| > \frac{\delta}{2} \right) \xrightarrow[N \to \infty]{} 0 \quad (5.3)$$

and

$$\mu_N\left(\left| \frac{1}{N^d} \sum_{x \in \mathbb{Z}^d} S_{0,tN^2}^N G(\tfrac{x}{N}) \eta(x) - \int_{\mathbb{R}^d} G(u) \rho_t^\Sigma(u) \, \mathrm{d} u \right| > \frac{\delta}{2} \right) \xrightarrow[N \to \infty]{} 0, \quad (5.4)$$

We do not prove tightness first because the computation done to prove (5.3) in Lemma 5.1 of Section 5.1 will be used again to prove tightness in Proposition 5.5 in Section 5.2.

Let us now shed more light on (5.3) and (5.4). Observe that the first term in the r.h.s. of (5.1) is deterministic – once $\eta_0 = \eta \in \{0, 1\}^{\mathbb{Z}^d}$ is fixed – whereas the second term has mean zero and contains all stochasticity derived from the stirring construction. Indeed, by (4.12), (4.5), (4.10) and (4.11), we get

$$\mathbb{E}_\eta[\mathsf{X}_t^N(G)] = \frac{1}{N^d} \sum_{x \in \mathbb{Z}^d} G(\tfrac{x}{N}) \mathbb{E}_\eta[\eta_{tN^2}(x)] = \frac{1}{N^d} \sum_{x \in \mathbb{Z}^d} G(\tfrac{x}{N}) \mathbb{E}_\eta[\eta(\widehat{X}_{0,tN^2}^x)]$$

$$= \frac{1}{N^d} \sum_{x \in \mathbb{Z}^d} G(\tfrac{x}{N}) \widehat{S}_{0,tN^2} \eta(x) = \frac{1}{N^d} \sum_{x \in \mathbb{Z}^d} [S^N_{0,tN^2} G(\tfrac{x}{N})] \eta(x).$$

Thus, the decomposition (5.1) can be written as

$$\mathsf{X}^N_t(G) = \mathbb{E}_\eta[\mathsf{X}^N_t(G)] + \left(\mathsf{X}^N_t(G) - \mathbb{E}_\eta[\mathsf{X}^N_t(G)]\right),$$

where the first term on the r.h.s. is the expectation of the empirical density field and the second one is "noise", i.e. the (stochastic) deviation from the mean; hence it satisfies

$$\mathbb{E}_\eta \left[ \frac{1}{N^d} \sum_{x \in \mathbb{Z}^d} \int_0^{tN^2} [S^N_{r,tN^2} G(\tfrac{x}{N})] \, \mathrm{d}M_r(\eta_{r^-}, x) \right] = 0.$$

Therefore, when deriving the hydrodynamic limit – basically a Weak Law of Large Numbers (WLLN) – the proof of (3.15) reduces to proving that, firstly, the "noise" vanishes in probability and, secondly, that the expectation – when initialized according to $\mu_N$ – converges to the correct deterministic limit corresponding to the macroscopic equation; that is (5.3) and (5.4), respectively.

## 5.1 Convergence of finite dimensional distributions

In the present section, we prove (3.15) by means of (5.3) and (5.4).

**Proof of (5.3).** The convergence (5.3) is a consequence of Chebyshev's inequality and the following lemma, derived through an adaptation of the proof of Lemma 12 in [41].

**Lemma 5.1.** *For all initial configurations $\eta \in \{0,1\}^{\mathbb{Z}^d}$ and $t \geq 0$, we have*

$$\mathbb{E}_\eta \left[ \left( \frac{1}{N^d} \sum_{x \in \mathbb{Z}^d} \int_0^{tN^2} S^N_{r,tN^2} G(\tfrac{x}{N}) \, \mathrm{d}M_r(\eta_{r^-}, x) \right)^2 \right] \xrightarrow[N \to \infty]{} 0.$$

*Proof.* By (4.16), we get

$$\frac{1}{N^d} \sum_{x \in \mathbb{Z}^d} \int_0^{tN^2} S^N_{r,tN^2} G(\tfrac{x}{N}) \, \mathrm{d}M_r(\eta_{r^-}, x)$$

$$= \frac{1}{N^d} \sum_{\{x,y\} \in E_d} \int_0^{tN^2} \left( S^N_{r,tN^2} G(\tfrac{x}{N}) - S^N_{r,tN^2} G(\tfrac{y}{N}) \right) (\eta_{r^-}(y) - \eta_{r^-}(x)) \, \mathrm{d}\overline{\mathcal{N}}_r(\{x,y\}).$$

Recall that the compensated Poisson processes $\{\overline{\mathcal{N}}_\cdot(\{x,y\}), \{x,y\} \in E_d\}$ defined in (4.2) are of bounded variation in view of Assumption 2.1 and, moreover, they are independent over bonds by the same property of the Poisson processes defined in (4.1). Thus by Itô's isometry for jump processes and the independence over the bonds of the Poisson processes in (4.1), we obtain

$$\mathcal{V}_{t,\eta}^N(G) := \mathbb{E}_\eta\left[\left(\frac{1}{N^d}\sum_{x\in\mathbb{Z}^d}\int_0^{tN^2} S_{r,tN^2}^N G(\tfrac{x}{N})\,\mathrm{d}M_r(\eta_{r^-},x)\right)^2\right]$$

$$= \frac{1}{N^{2d}}\sum_{\{x,y\}\in E_d}\int_0^{tN^2}\left(S_{r,tN^2}^N G(\tfrac{x}{N}) - S_{r,tN^2}^N G(\tfrac{y}{N})\right)^2 \xi_{r,\eta}(\{x,y\})\,\lambda_r(\{x,y\})\,\mathrm{d}r\,,$$

where $\xi_{r,\eta}(\{x,y\}) := \mathbb{E}_\eta\left[(\eta_{r^-}(y) - \eta_{r^-}(x))^2\right]$. Note that, for all $r \geq 0$ and $\eta \in \{0,1\}^{\mathbb{Z}^d}$, $\xi_{r,\eta}(\{x,y\}) \in [0,1]$. Then, firstly recall the definition of the random walk generator in (3.8) as well as (5.2), so that

$$\mathcal{V}_{t,\eta}^N(G) \leq \frac{1}{N^{2d}}\sum_{\{x,y\}\in E_d}\int_0^{tN^2}\left(S_{r,tN^2}^N G(\tfrac{x}{N}) - S_{r,tN^2}^N G(\tfrac{y}{N})\right)^2 \lambda_r(\{x,y\})\,\mathrm{d}r$$

$$= \frac{1}{N^{2d}}\sum_{x\in\mathbb{Z}^d}\int_0^{tN^2} S_{r,tN^2}^N G(\tfrac{x}{N})(-A_r S_{r,tN^2}^N G)(\tfrac{x}{N})(x)\,\mathrm{d}r\,.$$

Secondly use Kolmogorov backward equation (A.6) for the forward transition semigroup, to obtain

$$\mathcal{V}_{t,\eta}^N(G) \leq \frac{1}{N^{2d}}\sum_{x\in\mathbb{Z}^d}\int_0^{tN^2} S_{r,tN^2}^N G(\tfrac{x}{N})\,\partial_r S_{r,tN^2}^N G(\tfrac{x}{N})\,\mathrm{d}r$$

$$= \frac{1}{N^{2d}}\sum_{x\in\mathbb{Z}^d}\int_0^{tN^2} \frac{1}{2}\partial_r\left(S_{r,tN^2}^N G(\tfrac{x}{N})\right)^2\,\mathrm{d}r\,.$$

After integrating and using Proposition A.2(2), we further write

$$\frac{1}{N^{2d}}\sum_{x\in\mathbb{Z}^d}\int_0^{tN^2}\frac{1}{2}\partial_r\left(S_{r,tN^2}^N G(\tfrac{x}{N})\right)^2\,\mathrm{d}r$$

$$= \frac{1}{2N^d}\cdot\frac{1}{N^d}\sum_{x\in\mathbb{Z}^d}\left(G(\tfrac{x}{N})^2 - \left(S_{0,tN^2}^N G(\tfrac{x}{N})\right)^2\right)$$

$$\leq \frac{1}{2N^d}\cdot\frac{1}{N^d}\sum_{x\in\mathbb{Z}^d} G(\tfrac{x}{N})^2\,. \tag{5.5}$$

Because $\frac{1}{N^d}\sum_{x\in\mathbb{Z}^d}G(\tfrac{x}{N})^2 \to \int_{\mathbb{R}^d} G(u)^2\,\mathrm{d}u < \infty$ as $N\to\infty$, and since $\mathcal{V}_{t,\eta}^N(G) \geq 0$, the conclusion follows. □

**Proof of (5.4).** Note that for proving (5.3) neither assumptions (a) nor (b) of Theorem 3.2 have been invoked. In what follows, the invariance principle of forward random walks, i.e. assumption (b) of Theorem 3.2, will play a crucial role. More precisely, we exploit conditions, given in terms of convergence of semigroups, that are equivalent to the invariance principle. In the time-homogeneous context, the correspondence between weak convergence of Feller processes and convergence of Feller semigroups is due to Trotter and Kurtz (see e.g. [15], [35]). For the sake of completeness, in the next theorem we point out how this correspondence translates in the time-inhomogeneous setting. In what follows, $C_0(\mathbb{R}^d)$ denotes the Banach space of real-valued continuous functions on $\mathbb{R}^d$ vanishing at infinity endowed with the sup norm $\|\cdot\|_\infty$.

**Theorem 5.2** (INVARIANCE PRINCIPLE). *The following statements are equivalent:*

(A) Weak convergence in path-space. *The forward random walks $\{X_{0,t}^x,\ x \in \mathbb{Z}^d,\ t \in [0,\infty)\}$ satisfy an* invariance principle with arbitrary starting positions *with covariance matrix $\Sigma$ (that is, assumption (b) in Theorem 3.2).*

(B) Uniform convergence of transition semigroups. *For all $T > 0$ and $G \in C_0(\mathbb{R}^d)$,*

$$\sup_{0 \le s \le t \le T}\ \sup_{x \in \mathbb{Z}^d} \left| S^N_{sN^2,tN^2} G(\tfrac{x}{N}) - \mathcal{S}^{\Sigma}_{t-s} G(\tfrac{x}{N}) \right| \xrightarrow[N \to \infty]{} 0, \tag{5.6}$$

*where $\{S^N_{s,t},\ t \in [s,\infty)\}$ with $s \in [0,\infty)$ is the forward semigroup defined in (5.2) and $\{\mathcal{S}^{\Sigma}_t,\ t \in [0,\infty)\}$ is the Brownian motion semigroup, introduced before (3.4).*

*An analogous equivalence holds for the backward random walks when replacing $\{X_{0,t}^\cdot,\ t \in [0,\infty)\}$ and $S^N_{sN^2,tN^2}$ by $\{\widehat{X}_{0,t}^\cdot,\ s \in [0,t]\}$ and $\widehat{S}^N_{sN^2,tN^2}$, respectively.*

We do not provide a detailed proof of Theorem 5.2, but we just mention the main lines. Firstly, by Assumption 2.1, the random walks under consideration are Feller processes (see Appendix A.2). Secondly, by viewing $\{\mathcal{S}^{\Sigma}_t,\ t \ge 0\}$ as an operator semigroup on $C_0(\mathbb{R}^d)$, the Schwartz space $\mathcal{S}(\mathbb{R}^d)$, being a dense and $\mathcal{S}^{\Sigma}_t$-invariant (for all $t \ge 0$) subset of $C_0(\mathbb{R}^d)$, is a core for the associated infinitesimal generator $\mathcal{A}^{\Sigma}$. As a consequence, the idea is to conclude by means of [33, Theorem 19.25] (up to required adaptations as e.g. in [15], Theorem 6.1 in Chapter 1 and Corollary 8.7 in Chapter 4, because pre-limit and limit processes do not take values in the same state space), which applies to the time-homogeneous setting, only. Hence, we first consider the transition semigroup for the (time-homogeneous) space-time process $\{(X^x_{s,s+\cdot},s+\cdot),\ x \in \mathbb{Z}^d,\ s \ge 0\}$ defined in Appendix A.2, we apply [33, Theorem 19.25] in this time-homogeneous setting and, then, by considering only functions $\widetilde{G} \in \mathcal{S}(\mathbb{R}^d \times (-\infty,\infty))$ which do not depend on the time-variable within a compact interval of $(-\infty,\infty)$ and smoothly vanish outside of it, we obtain Theorem 5.2.

Having an invariance principle for the forward random walks in the environment $\lambda$ allows to replace the uniform convergence (w.r.t. $x \in \mathbb{Z}^d$) in (5.6) with convergence in mean (w.r.t. the counting measure). The more precise statement is the content of the following proposition. We note that a similar result with a similar proof appears already in [41, Proposition 14] and [17, p. 536]. However, our statement differs from those just mentioned because we include a uniform convergence over time as well as we require it to hold for functions in $\mathscr{S}(\mathbb{R}^d)$ rather than for functions in $C_{\text{comp}}(\mathbb{R}^d)$ – the space of real-valued continuous compactly supported functions on $\mathbb{R}^d$. For the sake of completeness, we present its proof below.

**Proposition 5.3.** *Keep the same notation as in Theorem 5.2. Assume that condition* (b) *in Theorem 3.2 holds true. Then, for all $T > 0$ and $G \in \mathscr{S}(\mathbb{R}^d)$, we have*

$$\sup_{0 \leq s \leq t \leq T} \frac{1}{N^d} \sum_{x \in \mathbb{Z}^d} \left| S^N_{sN^2, tN^2} G(\tfrac{x}{N}) - \mathcal{S}^{\Sigma}_{t-s} G(\tfrac{x}{N}) \right| \xrightarrow[N \to \infty]{} 0 \,. \tag{5.7}$$

*Proof.* Note first that by the assumption of Proposition 5.3, (B) of Theorem 5.2 holds. Moreover, if we split $G = G^+ - G^-$ into its positive and negative parts, we have $G^\pm \in C_0(\mathbb{R}^d) \cap L^1(\mathbb{R}^d)$ and that there exist functions $H^\pm \in \mathscr{S}(\mathbb{R}^d)$ such that

$$0 \leq G^\pm(u) \leq H^\pm(u), \quad u \in \mathbb{R}^d \,. \tag{5.8}$$

Indeed, we take, for instance

$$H^\pm(u) := \widetilde{G}^\pm * \varphi_\varepsilon(u) = \int_{B_\varepsilon(0)} \widetilde{G}^\pm(u - v)\, \varphi_\varepsilon(v)\, dv$$

where $\varphi_\varepsilon : \mathbb{R}^d \to \mathbb{R}$ is the canonical mollifier in $C_{\text{comp}}(\mathbb{R}^d)$ and with support in $B_\varepsilon(0)$ – in particular, non-negative and whose integral equals one – and

$$\widetilde{G}^\pm(u) := \sup_{v \in B_\varepsilon(0)} G^\pm(u - v) \,.$$

With this choice, we have $H^\pm \in C^\infty(\mathbb{R}^d)$ and, moreover, $H^\pm$ decays at infinity like $G^\pm$. Hence, $H^\pm \in \mathscr{S}(\mathbb{R}^d)$. As a consequence, there exist constants $C^\pm > 0$ such that

$$\sup_{0 \leq t \leq T} |\mathcal{S}^{\Sigma}_t G^\pm(u)| \leq \frac{C^\pm}{1 + |u|^{2d}}, \quad u \in \mathbb{R}^d \,. \tag{5.9}$$

This follows from the bounds (5.8), the fact that $\mathcal{S}^{\Sigma}_t$ acts as convolution with a non-degenerate Gaussian kernel and the use of Fourier transformation in $\mathscr{S}(\mathbb{R}^d)$. Moreover, because of the uniform continuity of $G^\pm$, we have

$$\sup_{0 \leq t \leq T} \sup_{|u-v| < \delta} \left| \mathcal{S}^{\Sigma}_t G^\pm(u) - \mathcal{S}^{\Sigma}_t G^\pm(v) \right| \leq \sup_{0 \leq t \leq T} \sup_{|u-v| < \delta} \left| G^\pm(u) - G^\pm(v) \right| \xrightarrow[\delta \to 0]{} 0$$

and, thus, by the integrability of $G^\pm$,

$$\sup_{0 \leq t \leq T} \left| \frac{1}{N^d} \sum_{x \in \mathbb{Z}^d} \mathcal{S}^{\Sigma}_t G^\pm(\tfrac{x}{N}) - \int_{\mathbb{R}^d} \mathcal{S}^{\Sigma}_t G^\pm(u)\, du \right| \xrightarrow[N \to \infty]{} 0 \,. \tag{5.10}$$

Let us now prove

$$\sup_{0\leq s\leq t\leq T} \frac{1}{N^d} \sum_{x\in\mathbb{Z}^d} \left| S^N_{sN^2,tN^2} G^{\pm}(\tfrac{x}{N}) - \mathcal{S}^{\Sigma}_{t-s} G^{\pm}(\tfrac{x}{N}) \right| \xrightarrow[N\to\infty]{} 0, \qquad (5.11)$$

from which (5.7) follows.

Because $|c| = c + 2\max\{-c, 0\}$ for all $c \in \mathbb{R}$, we have

$$\sup_{0\leq s\leq t\leq T} \frac{1}{N^d} \sum_{x\in\mathbb{Z}^d} \left| S^N_{sN^2,tN^2} G^{\pm}(\tfrac{x}{N}) - \mathcal{S}^{\Sigma}_{t-s} G^{\pm}(\tfrac{x}{N}) \right|$$

$$\leq \sup_{0\leq s\leq t\leq T} \frac{1}{N^d} \sum_{x\in\mathbb{Z}^d} \left( S^N_{sN^2,tN^2} G^{\pm}(\tfrac{x}{N}) - \mathcal{S}^{\Sigma}_{t-s} G^{\pm}(\tfrac{x}{N}) \right)$$

$$+ \sup_{0\leq s\leq t\leq T} \frac{2}{N^d} \sum_{x\in\mathbb{Z}^d} \max\left\{ \mathcal{S}^{\Sigma}_{t-s} G^{\pm}(\tfrac{x}{N}) - S^N_{sN^2,tN^2} G^{\pm}(\tfrac{x}{N}), 0 \right\} . \qquad (5.12)$$

By using, on the one side, (4.8), (4.9) and $\sum_{\frac{x}{N}\notin \mathcal{K}_{\varepsilon}} \widehat{p}_{sN^2,tN^2}(y,x) \leq 1$, while, on the other side, $\int_{\mathbb{R}^d} \mathcal{S}^{\Sigma}_t G^{\pm}(u)\,du = \int_{\mathbb{R}^d} G^{\pm}(u)\,du$ and (5.10), the first term in the r.h.s. of (5.12) vanishes as $N \to \infty$. Indeed,

$$\sup_{0\leq s\leq t\leq T} \left| \frac{1}{N^d} \sum_{x\in\mathbb{Z}^d} \left( S^N_{sN^2,tN^2} G^{\pm}(\tfrac{x}{N}) - \mathcal{S}^{\Sigma}_{t-s} G^{\pm}(\tfrac{x}{N}) \right) \right|$$

$$= \sup_{0\leq s\leq t\leq T} \left| \frac{1}{N^d} \sum_{y\in\mathbb{Z}^d} G^{\pm}(\tfrac{y}{N}) \sum_{x\in\mathbb{Z}^d} \widehat{p}_{sN^2,tN^2}(y,x) - \frac{1}{N^d} \sum_{x\in\mathbb{Z}^d} \mathcal{S}^{\Sigma}_{t-s} G^{\pm}(\tfrac{x}{N}) \right|$$

$$\leq \left| \frac{1}{N^d} \sum_{x\in\mathbb{Z}^d} G^{\pm}(\tfrac{x}{N}) - \int_{\mathbb{R}^d} G^{\pm}(u)\,du \right|$$

$$+ \sup_{0\leq s\leq t\leq T} \left| \frac{1}{N^d} \sum_{x\in\mathbb{Z}^d} \mathcal{S}^{\Sigma}_t G^{\pm}(\tfrac{x}{N}) - \int_{\mathbb{R}^d} \mathcal{S}^{\Sigma}_t G^{\pm}(u)\,du \right| .$$

Moreover, we have, for all $N \in \mathbb{N}$ and $x \in \mathbb{Z}^d$,

$$\sup_{0\leq s\leq t\leq T} \max\left\{ \mathcal{S}^{\Sigma}_{t-s} G^{\pm}(\tfrac{x}{N}) - S^N_{sN^2,tN^2} G^{\pm}(\tfrac{x}{N}), 0 \right\} \leq \sup_{0\leq t\leq T} \mathcal{S}^{\Sigma}_t G^{\pm}(\tfrac{x}{N}) . \qquad (5.13)$$

Therefore, by (5.13), (5.9) and (5.6), we conclude

$$\limsup_{N\to\infty} \sup_{0\leq s\leq t\leq T} \frac{2}{N^d} \sum_{x\in\mathbb{Z}^d} \max\left\{ \mathcal{S}^{\Sigma}_{t-s} G^{\pm}(\tfrac{x}{N}) - S^N_{sN^2,tN^2} G^{\pm}(\tfrac{x}{N}), 0 \right\}$$

$$\leq 2\, V_M \limsup_{N\to\infty} \sup_{0\leq s\leq t\leq T} \sup_{|\tfrac{x}{N}|\leq M} \left| \mathcal{S}^{\Sigma}_{t-s} G^{\pm}(\tfrac{x}{N}) - S^N_{sN^2,tN^2} G^{\pm}(\tfrac{x}{N}) \right|$$

$$+ \limsup_{N \to \infty} \frac{2}{N^d} \sum_{|\frac{x}{N}| > M} \frac{C^{\pm}}{1 + |\frac{x}{N}|^{2d}} \;=\; 2 \int_{|u| > M} \frac{C^{\pm}}{1 + |u|^{2d}} \, du \xrightarrow[M \to \infty]{} 0,$$

where

$$V_M := \sup_{N \in \mathbb{N}} \frac{\operatorname{card}\left\{x \in \mathbb{Z}^d : |\tfrac{x}{N}| \leq M\right\}}{N^d} < \infty.$$

□

We apply Proposition 5.3 and assumption (a) of Theorem 3.2 to prove (5.4) and conclude the characterization of the finite-dimensional distributions of the limiting density field, that is, (3.15).

Let $\{\rho_t^\Sigma, \, t \geq 0\}$ be the unique weak solution of the Cauchy problem as given in (3.4). Moreover, note that $\mathcal{S}_t^\Sigma \mathscr{S}(\mathbb{R}^d) \subset \mathscr{S}(\mathbb{R}^d)$ for all $t \geq 0$. Hence, for any family of probability measures $\{\mu_N, N \in \mathbb{N}\}$ associated to the density profile $\rho_\bullet$ (see (3.10) for the definition), we obtain

$$\mu_N \left( \left| \frac{1}{N^d} \sum_{x \in \mathbb{Z}^d} \mathcal{S}_t^\Sigma G(\tfrac{x}{N}) \eta(x) - \int_{\mathbb{R}^d} \mathcal{S}_t^\Sigma G(u) \rho_\bullet(u) \, du \right| > \delta \right) \xrightarrow[N \to \infty]{} 0, \qquad (5.14)$$

for all $t \geq 0$ and all $\delta > 0$. In turn, (5.4) comes as a consequence of (5.14) and the following lemma.

**Lemma 5.4.** *For all $t \geq 0$, all $G \in \mathscr{S}(\mathbb{R}^d)$ and for any sequence of probability measures $\{\tilde{\mu}_N, \, N \in \mathbb{N}\}$ in $\{0,1\}^{\mathbb{Z}^d}$, we have, for all $\delta > 0$,*

$$\tilde{\mu}_N \left( \left| \frac{1}{N^d} \sum_{x \in \mathbb{Z}^d} \left( S_{0,tN^2}^N G(\tfrac{x}{N}) - \mathcal{S}_t^\Sigma G(\tfrac{x}{N}) \right) \eta(x) \right| > \delta \right) \xrightarrow[N \to \infty]{} 0. \qquad (5.15)$$

*Proof.* Because $\eta(x) \leq 1$, we obtain

$$\left| \frac{1}{N^d} \sum_{x \in \mathbb{Z}^d} \left( S_{0,tN^2}^N G(\tfrac{x}{N}) - \mathcal{S}_t^\Sigma G(\tfrac{x}{N}) \right) \eta(x) \right| \leq \frac{1}{N^d} \sum_{x \in \mathbb{Z}^d} \left| S_{0,tN^2}^N G(\tfrac{x}{N}) - \mathcal{S}_t^\Sigma G(\tfrac{x}{N}) \right|.$$

Then we obtain (5.15) via Proposition 5.3 together with Markov's inequality. □

## 5.2 Tightness

In this section we prove tightness of the sequence of density fields $\{\mathsf{X}_\cdot^N, \, N \in \mathbb{N}\}$ in the Skorokhod space $\mathsf{D}([0, T], \mathscr{S}'(\mathbb{R}^d))$. Note that tightness of the distributions $\{\mathsf{X}_\cdot^N, \, N \in \mathbb{N}\}$ is implied by tightness of the density fields evaluated at all functions $G \in \mathscr{S}(\mathbb{R}^d)$ (see [39]). Hence, it suffices to discuss tightness of the sequence $\{\mathsf{X}_\cdot^N(G), \, N \in \mathbb{N}\}$ in $\mathsf{D}([0, T], \mathbb{R})$, for all $G \in \mathscr{S}(\mathbb{R}^d)$.

The criterion we use is given in Appendix B, Theorem B.4. Note that we cannot use Aldous-Rebolledo criterion (see e.g. [34]), which relies ultimately on Doob's maximal martingale inequality. Indeed, instead of decomposing the empirical density fields into a predictable term and a martingale term, we employed the mild solution representation (4.17) for which maximal inequalities for martingales do not apply. We postpone to Appendix B any precise statements and anticipate that in our case the proof boils down to prove the following.

**Proposition 5.5** (TIGHTNESS OF THE EMPIRICAL DENSITY FIELDS). *For any initial configuration $\eta \in \{0,1\}^{\mathbb{Z}^d}$ and all $G \in \mathscr{S}(\mathbb{R}^d)$, the real-valued stochastic processes $\{X_\cdot^N(G), \ N \in \mathbb{N}\}$ satisfy the following conditions:*

(T1) *For all $t \in [0,T]$,*
$$\lim_{m\to\infty} \limsup_{N\to\infty} \mathbb{P}\left(|X_t^N(G)| > m\right) = 0.$$

(T2) *For all $\varepsilon > 0$, there exist values $h_\varepsilon > 0$ and $N_\varepsilon \in \mathbb{N}$ such that for all $N \geq N_\varepsilon$ we find deterministic functions $\psi_\varepsilon^N, \psi_\varepsilon : [0, h_\varepsilon] \to [0,\infty)$ and non-negative values $\phi_\varepsilon^N$ satisfying, for all $N \geq N_\varepsilon$, the following properties:*

  (i) *The functions $\psi_\varepsilon^N$ are non-decreasing.*
  (ii) *For all $h \in [0, h_\varepsilon]$ and $t \in [0,T]$, we have*
$$\mathbb{P}\left(|X_{t+h}^N(G) - X_t^N(G)| > \varepsilon \,|\, \mathsf{F}_t^N\right) \leq \psi_\varepsilon^N(h), \quad \mathbb{P}\text{-a.s.},$$
  *where $\{\mathsf{F}_t^N, \ t \geq 0\}$ denotes the natural filtration associated to $\{X_t^N, \ t \geq 0\}$.*
  (iii) *For all $h \in [0, h_\varepsilon]$, we have $\psi_\varepsilon^N(h) \leq \psi_\varepsilon(h) + \phi_\varepsilon^N$.*
  (iv) *$\phi_\varepsilon^N \to 0$ as $N \to \infty$.*
  (v) *$\psi_\varepsilon(h) \to 0$ as $h \to 0$.*

*As a consequence of Theorems B.4 and B.2 below and [39, Theorem 4.1], $\{X_\cdot^N, \ N \in \mathbb{N}\}$ is a tight sequence in $D([0,T], \mathscr{S}'(\mathbb{R}^d))$.*

*Proof.* Statement (T1) is a direct consequence of (3.15), that we proved in Section 5.1. In what follows, we prove (T2).

For all $N \in \mathbb{N}$ and $t, t+h \in [0,T]$, writing $X_{t+h}^N(G)$ via (5.1), then using for these terms (4.11), Chapman-Kolmogorov equation for $\{S_{s,t}, \ t \in [s, \infty)\}$ (see Proposition A.2(6)), and (4.17), we get the decomposition

$$X_{t+h}^N(G) - X_t^N(G) = \frac{1}{N^d} \sum_{x \in \mathbb{Z}^d} (S_{tN^2,(t+h)N^2}^N G(\tfrac{x}{N}) - G(\tfrac{x}{N})) \eta_{tN^2}(x)$$
$$+ \frac{1}{N^d} \sum_{x \in \mathbb{Z}^d} \int_{tN^2}^{(t+h)N^2} S_{r,(t+h)N^2}^N G(\tfrac{x}{N}) \, \mathrm{d}M_r(\eta_{r^-}, x).$$

Thus, we obtain

$$\mathbb{P}\left(|\mathsf{X}_{t+h}^N(G) - \mathsf{X}_t^N(G)| > \varepsilon \,\Big|\, \mathsf{F}_t^N\right) \leq \mathcal{X}_{t,t+h}^N(\varepsilon) + \mathcal{Y}_{t,t+h}^N(\varepsilon), \qquad (5.16)$$

where

$$\mathcal{X}_{t,t+h}^N(\varepsilon) := \mathbb{P}\left(\left|\frac{1}{N^d}\sum_{x\in\mathbb{Z}^d}(S_{tN^2,(t+h)N^2}^N G(\tfrac{x}{N}) - G(\tfrac{x}{N}))\eta_{tN^2}(x)\right| > \frac{\varepsilon}{2}\,\bigg|\,\mathsf{F}_t^N\right) \qquad (5.17)$$

$$\mathcal{Y}_{t,t+h}^N(\varepsilon) := \mathbb{P}\left(\left|\frac{1}{N^d}\sum_{x\in\mathbb{Z}^d}\int_{tN^2}^{(t+h)N^2} S_{r,(t+h)N^2}^N G(\tfrac{x}{N})\,\mathrm{d}M_r(\eta_{r^-},x)\right| > \frac{\varepsilon}{2}\,\bigg|\,\mathsf{F}_t^N\right) \qquad (5.18)$$

and we estimate separately the two terms $\mathcal{X}_{t,t+h}^N(\varepsilon)$ and $\mathcal{Y}_{t,t+h}^N(\varepsilon)$ in (5.17) and (5.18), respectively. We start with $\mathcal{X}_{t,t+h}^N(\varepsilon)$. The bound $\eta_t(x) \leq 1$ yields

$$\mathcal{X}_{t,t+h}^N(\varepsilon) \leq \mathbb{P}\left(\frac{1}{N^d}\sum_{x\in\mathbb{Z}^d}\left|S_{tN^2,(t+h)N^2}^N G(\tfrac{x}{N}) - G(\tfrac{x}{N})\right| > \frac{\varepsilon}{2}\,\bigg|\,\mathsf{F}_t^N\right)$$

and the probability on the r.h.s. vanishes as $N \to \infty$. This can be seen as follows:

($\alpha$) By Proposition 5.3, there exists a sufficiently large $N_\varepsilon \in \mathbb{N}$ such that, for all $N \geq N_\varepsilon$, we have

$$\sup_{0\leq t\leq t+h\leq T}\frac{1}{N^d}\sum_{x\in\mathbb{Z}^d}\left|S_{tN^2,(t+h)N^2}^N G(\tfrac{x}{N}) - \mathcal{S}_h^\Sigma G(\tfrac{x}{N})\right| \leq \frac{\varepsilon}{4}\,; \qquad (5.19)$$

($\beta$) By the strong continuity of $\{\mathcal{S}_h^\Sigma, h \geq 0\}$ and the uniform integrability of $\{\mathcal{S}_h^\Sigma G, h \in [0,T]\}$ also used in the proof of Proposition 5.3, one can show that there exists $h_\varepsilon > 0$ – independent of $N \in \mathbb{N}$ – such that, for all $h \in [0, h_\varepsilon]$ and $N \geq N_\varepsilon$, we have

$$\frac{1}{N^d}\sum_{x\in\mathbb{Z}^d}\left|\mathcal{S}_h^\Sigma G(\tfrac{x}{N}) - G(\tfrac{x}{N})\right| \leq \frac{\varepsilon}{4}\,. \qquad (5.20)$$

We then obtain from (5.19) and (5.20), for all $h \in [0, h_\varepsilon]$ and $N \geq N_\varepsilon$,

$$\frac{1}{N^d}\sum_{x\in\mathbb{Z}^d}\left|S_{tN^2,(t+h)N^2}^N G(\tfrac{x}{N}) - G(\tfrac{x}{N})\right| \leq \frac{\varepsilon}{2}\,.$$

Hence, for all $N \geq N_\varepsilon$ and $h \in [0, h_\varepsilon]$,

$$\mathcal{X}_{t,t+h}^N(\varepsilon) = 0\,. \qquad (5.21)$$

To bound $\mathcal{Y}_{t,t+h}^N(\varepsilon)$ in (5.18), we combine Chebyshev's inequality and the argument in the proof of Lemma 5.1 (which gave (5.5) – note that we applied Itō's isometry for the conditional expectation) to get

$$\mathcal{Y}^N_{t,t+h}(\varepsilon) \leq \frac{4}{\varepsilon^2} \cdot \frac{1}{2N^d} \cdot \frac{1}{N^d} \sum_{x \in \mathbb{Z}^d} \left( G^2(\tfrac{x}{N}) - (S^N_{tN^2,(t+h)N^2} G)^2(\tfrac{x}{N}) \right) . \tag{5.22}$$

Recall the values $N_\varepsilon \in \mathbb{N}$ and $h_\varepsilon > 0$ obtained from conditions (5.19) and (5.20). For all $N \geq N_\varepsilon$ and $h \in [0, h_\varepsilon]$, define the function $\psi^N_\varepsilon : [0, h_\varepsilon] \to [0, \infty)$ as

$$\psi^N_\varepsilon(h) = \left( \frac{4\, C_G}{\varepsilon^2 N^d} \right) \mathcal{Z}^N_h , \tag{5.23}$$

where $C_G$ is given by

$$C_G := \sup_{N \in \mathbb{N}} \frac{1}{N^d} \sum_{x \in \mathbb{Z}^d} |G(\tfrac{x}{N})| < \infty , \tag{5.24}$$

and

$$\mathcal{Z}^N_h := \sup_{t \in [0,T]} \sup_{h' \in [0,h]} \sup_{x \in \mathbb{Z}^d} \left| G(\tfrac{x}{N}) - S^N_{tN^2,(t+h')N^2} G(\tfrac{x}{N}) \right| . \tag{5.25}$$

We observe that, for all $N \geq N_\varepsilon$, $\psi^N_\varepsilon$ is non-decreasing: indeed $\mathcal{Z}^N_{h'} \leq \mathcal{Z}^N_{h''}$ if $h' \leq h''$, given that $h', h'' \in [0, h_\varepsilon]$. This yields (i) for the functions $\{\psi^N_\varepsilon, N \geq N_\varepsilon\}$.

Now we prove (ii). We claim that $\mathcal{Y}^N_{t,t+h}(\varepsilon) \leq \psi^N_\varepsilon(h)$. Indeed, by (5.22),

$$\frac{4}{\varepsilon^2} \cdot \frac{1}{2N^d} \cdot \frac{1}{N^d} \sum_{x \in \mathbb{Z}^d} \left( G^2(\tfrac{x}{N}) - (S^N_{tN^2,(t+h)N^2} G)^2(\tfrac{x}{N}) \right)$$

$$\leq \frac{4}{\varepsilon^2} \cdot \frac{1}{2N^d} \cdot \frac{1}{N^d} \sum_{x \in \mathbb{Z}^d} |G(\tfrac{x}{N})| \cdot |G(\tfrac{x}{N}) - S^N_{tN^2,(t+h)N^2} G(\tfrac{x}{N})|$$

$$+ \frac{4}{\varepsilon^2} \cdot \frac{1}{2N^d} \cdot \frac{1}{N^d} \sum_{x \in \mathbb{Z}^d} |S^N_{tN^2,(t+h)N^2} G(\tfrac{x}{N})| \cdot |G(\tfrac{x}{N}) - S^N_{tN^2,(t+h)N^2} G(\tfrac{x}{N})|$$

$$\leq \left( \frac{4}{\varepsilon^2} \cdot \frac{1}{2N^d} \cdot \frac{1}{N^d} \sum_{x \in \mathbb{Z}^d} \left( |G(\tfrac{x}{N})| + |S^N_{tN^2,(t+h)N^2} G(\tfrac{x}{N})| \right) \right) \cdot \mathcal{Z}^N_h \leq \psi^N_\varepsilon(h), \tag{5.26}$$

where in the last inequality we used that, by (4.9) and (4.8), we have

$$\sum_{x \in \mathbb{Z}^d} |S^N_{tN^2,(t+h)N^2} G(\tfrac{x}{N})| \leq \sum_{x \in \mathbb{Z}^d} S^N_{tN^2,(t+h)N^2} |G|(\tfrac{x}{N}) = \sum_{x \in \mathbb{Z}^d} |G(\tfrac{x}{N})| \leq N^d C_G .$$

As a consequence, for our choices of $N_\varepsilon$ and $h_\varepsilon$ (see also (5.16) and its bounds, (5.21) and (5.26)), we have

$$\mathbb{P}\left( |\mathsf{X}^N_{t+h}(G) - \mathsf{X}^N_t(G)| > \varepsilon \,\big|\, \mathsf{F}^N_t \right) \leq \mathcal{X}^N_{t,t+h}(\varepsilon) + \mathcal{Y}^N_{t,t+h}(\varepsilon) \leq \psi^N_\varepsilon(h),$$

and, in turn, (ii).

Now we prove the last three items, namely (iii), (iv) and (v). By the triangle inequality, we obtain

$$\mathcal{Z}_h^N \leq \sup_{h' \in [0,h]} \sup_{u \in \mathbb{R}^d} |G(u) - \mathcal{S}_{h'}^\Sigma G(u)| + \sup_{0 \leq s \leq t \leq T} \sup_{x \in \mathbb{Z}^d} |S_{sN^2,tN^2}^N G(\tfrac{x}{N}) - \mathcal{S}_{t-s}^\Sigma G(\tfrac{x}{N})| \,. \tag{5.27}$$

This leads us to the following definitions: for all $h \in [0, h_\varepsilon]$ and $N \geq N_\varepsilon$,

$$\psi_\varepsilon(h) := \left(\frac{4\,C_G}{\varepsilon^2 (N_\varepsilon)^d}\right) \sup_{h' \in [0,h]} \sup_{u \in \mathbb{R}^d} |G(u) - \mathcal{S}_{h'}^\Sigma G(u)| \tag{5.28}$$

$$\phi_\varepsilon^N := \left(\frac{4\,C_G}{\varepsilon^2 N^d}\right) \sup_{0 \leq s \leq t \leq T} \sup_{x \in \mathbb{Z}^d} |S_{sN^2,tN^2}^N G(\tfrac{x}{N}) - \mathcal{S}_{t-s}^\Sigma G(\tfrac{x}{N})| \,. \tag{5.29}$$

As a consequence, (iii) holds.

We obtain (iv), i.e. $\phi_\varepsilon^N \to 0$ as $N \to \infty$, from (5.6), i.e. forward semigroup uniform convergence. Alternatively, the contraction property of the semigroups $\{S_{s,t},\ t \in [s, \infty)\}$ and $\{\mathcal{S}_t^\Sigma,\ t \geq 0\}$ (cf. Proposition A.2(4)), yields

$$\sup_{0 \leq s \leq t \leq T} \sup_{x \in \mathbb{Z}^d} |S_{sN^2,tN^2}^N G(\tfrac{x}{N}) - \mathcal{S}_{t-s}^\Sigma G(\tfrac{x}{N})| \leq 2 \sup_{u \in \mathbb{R}^d} |G(u)| \,. \tag{5.30}$$

By combining (5.30) with $\phi_\varepsilon^N \geq 0$ and $\phi_\varepsilon^N \leq (\sup_u |G(u)| \tfrac{8}{\varepsilon^2} C_G) \tfrac{1}{N^d}$ leads to (iv).

The property $\psi_\varepsilon(h) \to 0$ as $h \to 0$, i.e. item (v), is a consequence of the strong continuity of the contraction semigroup $\{\mathcal{S}_t^\Sigma,\ t \geq 0\}$. This concludes the proof. □

## A  Time-inhomogeneous random walks: graphical construction and properties

In this appendix we collect some basic facts about time-inhomogeneous random walks. In particular, first we detail a dynamic version of Harris graphical construction ([29]) based on a percolation argument, which was summarized in Section 4.1. Then, we prove Proposition 4.3. In conclusion, we show that the random walks obtained are indeed Feller processes. We rely on the notation in Section 4.1.

### A.1  Graphical construction of random walks

In this section we explain in detail the graphical construction which defines the percolation structure on which we build the families of forward and backward random walks given in Section 4.1. The main difficulty comes from the loss of space-time translation invariance due to the dynamic environment $\lambda$. We deal with this difficulty by using Assumption 2.1 about the uniform boundedness of the conductances by $\mathfrak{a} > 0$. It will enable us to relate the percolation structure built from the inhomogeneous Poisson processes to a *bond percolation model in $\mathbb{Z}^d$* ([28]). Using the latter, we can construct the families of random walks by piecing together paths defined on sufficiently small time intervals which cover the whole positive real line.

**Remark A.1.** *The uniform boundedness assumption could in principle be relaxed as long as results from bond percolation models transfer to our inhomogeneous setting. For examples of weaker assumptions on the conductances when $d = 1$, see e.g. [17, Lemma 2.1] or [6].*

**Stochastic domination.** Let $\{\mathfrak{R}.(\{x,y\}), \{x,y\} \in E_d\}$ be a family of i.i.d. Poisson processes of intensity $\mathfrak{a}$ defined on the probability space $(\Xi, \mathfrak{F}, \{\mathfrak{F}_t, t \geq 0\}, \mathfrak{P})$, where $\{\mathfrak{F}_t, t \geq 0\}$ denotes its natural filtration and $\mathfrak{F} = \sigma(\cup_{t \geq 0} \mathfrak{F}_t)$. By a thinning procedure (see e.g. [37]), we construct the family of inhomogeneous Poisson processes $\{\mathcal{N}.(\{x,y\}), \{x,y\} \in E_d\}$ given in (4.1) as follows: for all $n \in \mathbb{N}$ and $\{x,y\} \in E_d$, if we denote by $T_n(\{x,y\})$ the random time at which the $n$-th event of $\mathfrak{R}.(\{x,y\})$ has occurred, we erase this random time with probability

$$1 - \mathfrak{a}^{-1} \cdot \lambda_{T_n(\{x,y\})}(\{x,y\}) .$$

We proceed analogously and independently for all random times $\{T_m(\{x,y\}), m \in \mathbb{N}, \{x,y\} \in E_d\}$. We denote the probability space induced by $\{\mathfrak{R}.(\{x,y\}), \{x,y\} \in E_d\}$ and this thinning procedure by $(\Omega, \mathcal{F}, \{\mathcal{F}_t, t \geq 0\}, \mathbb{P})$, where $\{\mathcal{F}_t, t \geq 0\}$ denotes the corresponding natural filtration and $\mathcal{F} = \sigma\left(\bigcup_{t \geq 0} \mathcal{F}_t\right)$. Then the remaining random points form the family of inhomogeneous Poisson process $\{\mathcal{N}.(\{x,y\}), \{x,y\} \in E_d\}$ introduced in (4.1), see also [37].

Given this construction, for all $\{x,y\} \in E_d$ and $t \geq s$, the number of Poissonian events of $\mathfrak{R}.(\{x,y\})$ in the time interval $[s,t]$ $\mathbb{P}$-a.s. dominates the number of events of $\mathcal{N}.(\{x,y\})$ in the same time interval.

**Percolation and active islands.** Let us first consider the family of i.i.d. Poisson processes $\{\mathfrak{R}.(\{x,y\}), \{x,y\} \in E_d\}$. For all $t \geq s$, we say that the bond $\{x,y\} \in E_d$ is *open in* $[s,t]$ if

$$\mathfrak{R}_t(\{x,y\}) - \mathfrak{R}_{s^-}(\{x,y\}) \geq 1 .$$

We call the connected components of the subgraph consisting of sites of $\mathbb{Z}^d$ and bonds that are open in $[s,t]$ *open clusters in* $[s,t]$.

Bonds $\{x,y\} \in E_d$ are open in $[s,t]$ independently of each other with probability $\mathfrak{p}_{s,t}(\mathfrak{a}) = 1 - e^{\mathfrak{a}(t-s)}$. Hence, for all $t \geq s$, this induces a *bond percolation model in $\mathbb{Z}^d$ with density* $\mathfrak{p}_{s,t}(\mathfrak{a})$. As a consequence of the existence of a critical probability $\mathfrak{p}(d) \in (0,1]$ for bond percolation in $\mathbb{Z}^d$ (see [28, p. 13] for the case $d = 1$ and [28, Theorems 1.10–11] for the case $d \geq 2$), for any $d \geq 1$ there exists a value $h_c(d, \mathfrak{a}) > 0$ such that for all $s, t \in [0, \infty)$ with $0 < t - s < h_c(d, \mathfrak{a})$, the open clusters in $[s,t]$ are all finite $\mathfrak{P}$-almost surely. In particular, if we fix an interval size $\bar{h} < h_c(d, \mathfrak{a})$ and consider, for all $k \geq 0$, the interval $I_k = [k\bar{h}, (k+1)\bar{h}]$, we have

$$\mathfrak{P}(\text{for all } k \in \mathbb{N} \text{ all open clusters in } I_k \text{ are finite}) = 1 . \qquad (\text{A.1})$$

We turn now to the inhomogeneous Poisson processes $\{\mathcal{N}.(\{x,y\}),\ \{x,y\} \in E_d\}$. A bond $\{x,y\} \in E_d$ such that

$$\mathcal{N}_t(\{x,y\}) - \mathcal{N}_{s^-}(\{x,y\}) \geq 1$$

is said to be *active in* $[s,t]$. Note that, due to the thinning procedure, active bonds are open, but not necessarily vice versa. Hence, the connected components of the subgraph consisting of sites of $\mathbb{Z}^d$ and bonds that are active in $[s,t]$, denoted as *active island in* $[s,t]$, are $\mathbb{P}$-a.s. finite if $0 < t - s < h_c(d,\mathfrak{a})$. Combining this with (A.1), we get

$$\mathbb{P}(\text{for all } k \in \mathbb{N} \text{ all active islands in } I_k \text{ are finite}) = 1. \tag{A.2}$$

**Random walks.** Let us denote by $\mathcal{G}_{[s,t]}(x)$ the active island in $[s,t]$ containing the site $x \in \mathbb{Z}^d$. Then, by rephrasing (A.2), we have

$$\mathbb{P}(\text{for all } k \in \mathbb{N} \text{ and } x \in \mathbb{Z}^d, \mathcal{G}_{I_k}(x) \text{ is finite}) = 1. \tag{A.3}$$

Moreover, $\mathbb{P}$-a.s. each $\mathcal{G}_{I_k}(x)$ contains at most finitely-many Poissonian marks (and no marks at the times $k\bar{h}$, for all $k \in \mathbb{N}$).

As a consequence, for $\mathbb{P}$-a.e. $\omega \in \Omega$, if we choose $s, t \in I_k$ for some $k \in \mathbb{N}$ with $0 < t - s < h_c(d,\mathfrak{a})$, the random walks' paths $\{X^x_{s,r}[\omega], x \in \mathbb{Z}^d, r \in [s,t]\}$ and $\{\widehat{X}^y_{r,t}[\omega], y \in \mathbb{Z}^d, r \in [s,t]\}$ are all simultaneously well-defined. Indeed, for all $x$ and $y \in \mathbb{Z}^d$, it suffices to consider only finitely many Poissonian marks within $\mathcal{G}_{I_k}(x)[\omega]$ and $\mathcal{G}_{I_k}(y)[\omega]$, respectively, when performing the jumps (right-continuous for the forward random walks and left-continuous for the backward random walks). This procedure uniquely defines $X^x_{s,t}[\omega]$ and $\widehat{X}^y_{s,t}[\omega]$ simultaneously for all $x$ and $y \in \mathbb{Z}^d$ (see also Section 4).

If $s, t \in [0,\infty)$ with $t - s > 0$ belong to different intervals $I_{k(s)}$ and $I_{k(t)}$, respectively, with $k(s) < k(t)$, then, by piecing together the well-defined paths

$$\{X^x_{s,r}[\omega], r \geq s \text{ and } r \in I_{k(s)}\}, \ldots, \{X^x_{r,t}[\omega], r \leq t \text{ and } r \in I_{k(t)}\}$$

in ascending order w.r.t. $k \in \{k(s), \ldots, k(t)\}$ and

$$\{\widehat{X}^y_{r,t}[\omega], r \leq t \text{ and } r \in I_{k(t)}\}, \ldots, \{\widehat{X}^y_{s,r}[\omega], r \geq s \text{ and } r \in I_{k(s)}\}$$

in descending order w.r.t. $k \in \{k(s), \ldots, k(t)\}$, we obtain $X^x_{s,t}[\omega]$ and $\widehat{X}^y_{s,t}[\omega]$ for all $x$ and $y \in \mathbb{Z}^d$.

The property of the inhomogeneous Poisson processes $\{\mathcal{N}.(\{z,v\}), \{z,v\} \in E_d\}$ for which past and future are independent conditioned on the present state and our construction rules of the random walks imply that the processes $\{X^x_{s,t}, t \in [s,\infty)\}$ and $\{\widehat{X}^y_{s,t}, s \in [0,t]\}$ for all $x$ and $y \in \mathbb{Z}^d$ are Markovian w.r.t. the induced natural filtrations. This justifies the introduction in Section 4 of the transition probabilities $\{p_{s,t}(x,y), x,y \in \mathbb{Z}^d\}$ and $\{\widehat{p}_{s,t}(y,x), x,y \in \mathbb{Z}^d\}$, as well as of the semigroups

$\{S_{s,t},\ t \in [s,\infty)\}$ and $\{\widehat{S}_{s,t},\ s \in [0,t]\}$.

## A.2 Feller transition semigroups and generators

We study properties of the transition semigroups $\{S_{s,t},\ t \in [s,\infty)\}$ and $\{\widehat{S}_{s,t},\ s \in [0,t]\}$ introduced in (4.9) and (4.10) and their associated infinitesimal generators solving the associated Kolmogorov forward and backward equations as in (A.5) and (A.6), which turn out to be $\{A_t,\ t \geq 0\}$ and $\{A_{t^-},\ t \geq 0\}$ defined in (3.8) and (3.9). Indeed, for all $x, y \in \mathbb{Z}^d$ with the convention $\lambda_t(\{x,y\}) = 0$ if $\{x,y\} \notin E_d$, we have

$$\lim_{t \downarrow s} \frac{p_{s,t}(x,y)}{t-s} = \lambda_{s^+}(\{x,y\}) \quad \text{and} \quad \lim_{s \uparrow t} \frac{\widehat{p}_{s,t}(x,y)}{t-s} = \lambda_{t^-}(\{x,y\}), \qquad (A.4)$$

where the following limits $\lambda_{s^{\pm}}(\{x,y\}) = \lim_{h \downarrow 0} \lambda_{s \pm h}(\{x,y\})$ exist and, as the conductances are assumed to be càdlàg, $\lambda_{s^+}(\{x,y\}) = \lambda_s(\{x,y\})$.

In what follows, for a differentiable function $\phi : (-\infty, \infty) \to (X, \mathfrak{d})$, with $(X, \mathfrak{d})$ a metric space, we define

$$\partial_\tau \phi(\tau) = \lim_{h \downarrow 0} \tfrac{1}{h}(\phi(\tau+h) - \phi(\tau)) \quad \text{and} \quad \partial_{\tau^-} \phi(\tau) = \lim_{h \downarrow 0} \tfrac{1}{h}(\phi(\tau) - \phi(\tau-h)).$$

Moreover, $C_0(\mathbb{R}^d)$, resp. $C_0(\mathbb{Z}^d)$, denotes the Banach space of real-valued continuous functions on $\mathbb{R}^d$, resp. $\mathbb{Z}^d$, vanishing at infinity endowed with the sup norm $\|\cdot\|_\infty$.

The proofs of the next two propositions, which follow from Assumption 2.1, are left to the reader. For notational convenience, we extend the definitions of conductances, transition semigroups and generators to negative times.

**Proposition A.2** (TRANSITION SEMIGROUPS). *For all $f \in C_0(\mathbb{Z}^d)$ and $r \in [s,t]$, the following hold true:*

(1) Operators on $C_0(\mathbb{Z}^d)$. $S_{s,t}f \in C_0(\mathbb{Z}^d)$ and $\widehat{S}_{s,t}f \in C_0(\mathbb{Z}^d)$.

(2) Identity. $S_{t,t}f = \widehat{S}_{t,t}f = f$.

(3) Positivity. *If $f \geq 0$, then $S_{s,t}f \geq 0$ and $\widehat{S}_{s,t}f \geq 0$.*

(4) Contraction. $\|S_{s,t}f\|_\infty \leq \|f\|_\infty$ and $\|\widehat{S}_{s,t}f\|_\infty \leq \|f\|_\infty$.

(5) Conservativity. $S_{s,t}1 = 1$ and $\widehat{S}_{s,t}1 = 1$.

(6) Chapman-Kolmogorov equation.

$$S_{s,r}S_{r,t}f = S_{s,t}f \qquad \widehat{S}_{r,t}\widehat{S}_{s,r}f = \widehat{S}_{s,t}f.$$

(7) Strong continuity. *For all $T > 0$, $\lim_{h \downarrow 0} \sup_{0 \leq s \leq T} \|S_{s,s+h}f - f\|_\infty = 0$.*

**Proposition A.3** (INFINITESIMAL GENERATORS). *For all $f \in C_0(\mathbb{Z}^d)$ and $t \in (-\infty, \infty)$, the following hold true:*

(1) Domain. $A_t f$ and $A_{t^-} f \in C_0(\mathbb{Z}^d)$.

(2) Kolmogorov forward and backward equations.

$$\begin{aligned} \partial_t S_{s,t} f &= S_{s,t} A_t f & \partial_t \widehat{S}_{s,t} f &= A_t \widehat{S}_{s,t} f \\ \partial_{t^-} S_{s,t} f &= S_{s,t} A_{t^-} f & \partial_{t^-} \widehat{S}_{s,t} f &= A_{t^-} \widehat{S}_{s,t} f \end{aligned} \qquad (A.5)$$

*and*

$$\begin{aligned} \partial_s S_{s,t} f &= -A_s S_{s,t} f & \partial_s \widehat{S}_{s,t} f &= -\widehat{S}_{s,t} A_s f \\ \partial_{s^-} S_{s,t} f &= -A_{s^-} S_{s,t} f & \partial_{s^-} \widehat{S}_{s,t} f &= -\widehat{S}_{s,t} A_{s^-} f \,, \end{aligned} \qquad (A.6)$$

*where derivatives are meant w.r.t.* $\|\cdot\|_\infty$.

**Feller property.** We now consider the space-time process [51, Section 8.5.5]

$$\{(X^x_{s,s+\cdot}, s+\cdot),\ x \in \mathbb{Z}^d,\ s \in (-\infty, \infty)\} \qquad (A.7)$$

$$\{(\widehat{X}^y_{t-\cdot,t}, t-\cdot),\ y \in \mathbb{Z}^d,\ t \in (-\infty, \infty)\} \qquad (A.8)$$

associated to forward and backward random walks, respectively. These processes are time-homogeneous Markov processes on the state space $\mathbb{Z}^d \times (-\infty, \infty)$ with Markov pre-generators $\mathcal{B}$ and $\widehat{\mathcal{B}}$ given by

$$\mathcal{B} f(x,s) = A_s f(x,s) + \partial_s f(x,s) \qquad (A.9)$$

$$\widehat{\mathcal{B}} f(x,s) = A_{s^-} f(x,s) + \partial_s f(x,s)\,, \qquad (A.10)$$

for all

$$f \in \{f : f(x,\cdot) \in C^1_{\text{comp}}(-\infty, \infty) \text{ for all } x \in \mathbb{Z}^d\} \subset C_0(\mathbb{Z}^d \times (-\infty, \infty))$$

(see e.g. [8], [46, Chapter III.2], [51, Section 8.5.5]). Their closure in $C_0(\mathbb{Z}^d \times (-\infty, \infty))$ w.r.t. the supremum norm generates a strongly continuous contraction semigroup (see e.g. [14, Section II.3.28]).

Hence, by passing to this formulation, Propositions A.2 and A.3 guarantee that the forward and backward random walks are Feller processes, i.e. in the sense that their associated space-time processes are Feller processes as in [33, Chapter 19. Conditions (F1)–(F3)].

## A.3 Proof of Proposition 4.3

We have to show in this section that the infinite summation on the r.h.s. of (4.17) is (for $\mathbb{P}$-a.e. $\omega \in \Omega$) absolutely convergent and that it equals $\eta_t(x)[\omega]$. The proof relies on the construction of *active islands* (introduced in Appendix A.1) and on a finer control on their radius, which allows to obtain exponential bounds on the transition probabilities of the random walks. As a consequence, we prove identity (4.17) for all initial conditions $\eta \in \{0,1\}^{\mathbb{Z}^d}$ and all times $t \geq 0$.

The plan is the following. First we show that, when restricting to a finite summation, formula (4.17) indeed holds for $\mathbb{P}$-a.e. $\omega \in \Omega$. Then, based only on a percolation result on the radius of active islands in sufficiently small time intervals ([28]) and the uniform boundedness assumption of the conductances (Assumption 2.1), we obtain an exponential upper bound for the heat kernel. In conclusion, we prove that, for all initial conditions $\eta \in \{0,1\}^{\mathbb{Z}^d}$, for $\mathbb{P}$-a.e. $\omega \in \Omega$, the infinite summation in (4.17) is absolutely convergent, hence a rearrangement of the order of the summation, which does not change its value, gives us the result.

**Finite summations.** Among all *active islands in* $[s,t]$, namely the connected subgraphs of $(\mathbb{Z}^d, E_d)$ consisting of sites of $\mathbb{Z}^d$ and bonds $\{y,z\} \in E_d$ for which Poissonian events occurred in the time window $[s,t]$, i.e. for which

$$\mathcal{N}_t(\{y,z\}) - \mathcal{N}_{s^-}(\{y,z\}) \geq 1,$$

we denote by $\mathcal{G}_{[s,t]}[x]$ the unique active island in $[s,t]$ containing $x \in \mathbb{Z}^d$. Due to Assumption 2.1 (see Appendix A.1 for the detailed argument), there exists $h_c(d,\mathfrak{a}) > 0$ such that, for $\mathbb{P}$-a.e. realization $\omega \in \Omega$, for all $x \in \mathbb{Z}^d$ and all $s,t \in [0,\infty)$ with $0 < t - s < h_c(d,\mathfrak{a})$, the active island $\mathcal{G}_{[s,t]}(x)[\omega]$ is finite. As a consequence, both trajectories $\{X^x_{s,r}[\omega],\ r \in [s,t]\}$ and $\{\widehat{X}^x_{r,t}[\omega],\ r \in [s,t]\}$ are well-defined (see Section 4.1 and Appendix A.1). For the same reason, given $\{\eta_s(x)[\omega],\ x \in \mathbb{Z}^d\}$, the definition (4.12), i.e. $\eta_t(x)[\omega] = \eta_s(\widehat{X}^x_{s,t}[\omega])[\omega]$ in terms of the stirring process, poses no problem.

In the following lemma, due to the finiteness of active islands, we can give a precise meaning to (4.17) when restricting the summation only to particle positions within the same active island.

**Lemma A.4.** *Fix $x, z \in \mathbb{Z}^d$ and $s, t \in [0, \infty)$ with $0 < t - s < h_c(d,\mathfrak{a})$. Then, for $\mathbb{P}$-a.e. $\omega \in \Omega$ and any configuration $\eta_s \in \{0,1\}^{\mathbb{Z}^d}$, we have*

$$\sum_{y \in \mathcal{G}_{[s,t]}(z)[\omega]} \left( \widehat{p}_{s,t}(x,y)\eta_s(y) + \int_s^t \widehat{p}_{r,t}(x,y)\,dM_r(\eta_{r^-}[\omega],y)[\omega] \right)$$
$$= \begin{cases} \eta_t(x)[\omega] & \text{if } x \in \mathcal{G}_{[s,t]}(z)[\omega] \\ 0 & \text{otherwise}, \end{cases} \quad (\text{A.11})$$

*where $\eta_t(x)[\omega] = \eta_s(\widehat{X}^x_{s,t}[\omega])$.*

*Proof.* For notational convenience, let us set $s = 0$ and $t < h_c(d,\mathfrak{a})$. By recalling the definition of $dM_r$ in (4.16) and the following backward master equation (obtained by using (A.6))

$$\partial_r \widehat{p}_{r,t}(x,\cdot)(y) = -\sum_{v:\{y,v\} \in E_d} \lambda_r(\{y,v\})(\widehat{p}_{r,t}(x,v) - \widehat{p}_{r,t}(x,y)), \quad (\text{A.12})$$

we rearrange the l.h.s. in (A.11) to obtain:

$$\sum_{y \in \mathcal{G}_{[0,t]}(z)[\omega]} \left( \widehat{p}_{0,t}(x,y)\eta(y) + \int_0^t \eta_{r^-}(y)[\omega]\, \partial_r \widehat{p}_{r,t}(x,\cdot)(y)\, dr \right.$$

$$\left. + \int_0^t \eta_{r^-}(y)[\omega] \sum_{v:\{y,v\} \in E_d} (\widehat{p}_{r,t}(x,v) - \widehat{p}_{r,t}(x,y))\, d\mathcal{N}_r(\{y,v\})[\omega] \right). \quad \text{(A.13)}$$

Now, for all $y \in \mathcal{G}_{[0,t]}(z)[\omega]$, we denote by

$$0 \le s_1(y)[\omega] < \ldots < s_{n(y)[\omega]}(y)[\omega] \le t$$

the $n(y)[\omega]$ jump times occurred in a bond incident to $y$ in the time interval $[0,t]$, with the convention

$$s_0 = 0 \quad \text{and} \quad s_{n(y)[\omega]+1} = t.$$

Note that Assumption 2.1 assures that $n(y)[\omega] < \infty$. By recalling definition (4.12), each $y$-term in (A.13) admits the following decomposition (we write $s_k$ instead of $s_k(y)[\omega]$ for readability):

$$\widehat{p}_{0,t}(x,y)\eta(y) + \sum_{k=0}^{n(y)[\omega]} \eta(\widehat{X}^y_{0,s_k}[\omega]) \left(\widehat{p}_{s_{k+1},t}(x,y) - \widehat{p}_{s_k,t}(x,y)\right)$$

$$+ \sum_{k=0}^{n(y)[\omega]} \eta(\widehat{X}^y_{0,s_k}[\omega]) \left(\widehat{p}_{s_{k+1},t}(x, X^y_{s_k,s_{k+1}}[\omega]) - \widehat{p}_{s_{k+1},t}(x,y)\right),$$

which further simplifies as

$$\widehat{p}_{t,t}(x,y)\eta(\widehat{X}^y_{0,t}[\omega]) + \sum_{k=0}^{n(y)[\omega]-1} \left( \eta(\widehat{X}^y_{0,s_k}[\omega])\widehat{p}_{s_{k+1},t}(x, X^y_{s_k,s_{k+1}}[\omega]) \right.$$

$$\left. -\eta(\widehat{X}^y_{0,s_{k+1}}[\omega])\widehat{p}_{s_{k+1},t}(x,y) \right). \quad \text{(A.14)}$$

Now, for all $k = 0, \ldots, n(y)[\omega] - 1$, there exists a unique neighbor of $y$, here denoted by $v \in \mathbb{Z}^d$, for which $d\mathcal{N}_{s_{k+1}(y)}(\{y,v\})[\omega] = 1$. Note that $v \in \mathcal{G}_{[0,t]}(z)[\omega]$. As a consequence of the construction of the forward and backward random walks, we have

$$X^y_{s_k,s_{k+1}}[\omega] = v, \quad \widehat{X}^y_{0,s_k}[\omega] = \widehat{X}^v_{0,s_{k+1}}[\omega] \quad \text{and} \quad \widehat{X}^v_{0,s_k}[\omega] = \widehat{X}^y_{0,s_{k+1}}[\omega].$$

In turn, there will be exactly one term in the following sum

$$\sum_{\ell=0}^{n(v)[\omega]-1} \left( \eta(\widehat{X}^v_{0,s_\ell}[\omega])\widehat{p}_{s_{\ell+1},t}(x, X^v_{s_\ell,s_{\ell+1}}[\omega]) - \eta(\widehat{X}^v_{0,s_{\ell+1}}[\omega])\widehat{p}_{s_{\ell+1},t}(x,v) \right)$$

which cancels the corresponding $k$-th term in (A.14). Hence, after reordering these finite summations, (A.13) reduces to the following

$$\sum_{y \in \mathcal{G}_{[0,t]}(z)[\omega]} \widehat{p}_{t,t}(x,y) \eta(\widehat{X}^{y}_{0,t}[\omega]) \, .$$

The observation that $\widehat{p}_{t,t}(x,y) = \mathbf{1}_{\{x=y\}}$ concludes the proof. □

**Radius of active islands and absolute convergence.** We start by presenting a key estimate, direct consequence of [28, Theorem 3.4] and Assumption 2.1, on the radius of active islands:

**Fact A.1** ([28, THEOREM 3.4]). *For all $s, t \in [0, \infty)$, with $0 < t - s < h_c(d, \mathfrak{a})$, there exists $\chi(t - s) > 0$, such that*

$$\mathbb{P}\left(\exists y \in \mathbb{Z}^d : |y - x| = n \text{ and } y \in \mathcal{G}_{[s,t]}(x)\right) \leq e^{-\chi(t-s)n} \,, \tag{A.15}$$

*for all $n \in \mathbb{N}$ and $x \in \mathbb{Z}^d$. In words, the probability that the active island in $(s, t]$ containing $x \in \mathbb{Z}^d$ contains at least one site at distance $n$ from $x \in \mathbb{Z}^d$ is smaller than $e^{-\chi(t-s)n}$, for all $n \in \mathbb{N}$. The function $\chi : (0, h_c(d, \mathfrak{a})) \to (0, \infty)$ can be chosen to be non-increasing.*

For all $x \in \mathbb{Z}^d$, $t \geq 0$ and $\eta \in \{0, 1\}^{\mathbb{Z}^d}$, we need to give a precise meaning to the infinite sum in (4.17) for $\mathbb{P}$-a.e. realization $\omega \in \Omega$. More precisely, we need to ensure that this infinite sum is absolutely convergent, allowing us to reorder the summation so as to sum over finite active islands (over space and time) first, and then, to apply Lemma A.4. This is the content of the following lemma.

**Lemma A.5.** *Fix $x \in \mathbb{Z}^d$, $t > 0$ and a partition $\{0 = t_0, t_1, \ldots, t_n, t_{n+1} = t\}$ of $[0, t]$ finer than $h_c(d, \mathfrak{a})$, i.e. $t_{k+1} > t_k$ and $t_{k+1} - t_k < h_c(d, \mathfrak{a})$, for $k = 0, \ldots, n$. Then:*

(1) *There exist two constants $C, \overline{\chi} > 0$ (depending only on $t > 0$ and the partition $\{t_0, \ldots, t_{n+1}\}$ of $[0, t]$) such that, for all $m \in \mathbb{N}$,*

$$\sum_{y : |y-x|=m} \sum_{z_1 \in \mathbb{Z}^d} \cdots \sum_{z_n \in \mathbb{Z}^d} p_{0,t_1}(x, z_1) \cdots p_{t_n,t}(z_n, y) \leq C \, e^{-\overline{\chi} m} \,. \tag{A.16}$$

*An analogous result (with the same constants) holds for the backward transition probabilities:*

$$\sum_{y : |y-x|=m} \sum_{z_n \in \mathbb{Z}^d} \cdots \sum_{z_1 \in \mathbb{Z}^d} \widehat{p}_{t_n,t}(x, z_n) \cdots \widehat{p}_{0,t_1}(z_1, y) \leq C \, e^{-\overline{\chi} m} \,. \tag{A.17}$$

(2) *For $\mathbb{P}$-a.e. $\omega \in \Omega$ and all initial configurations $\eta \in \{0, 1\}^{\mathbb{Z}^d}$,*

$$\sum_{z_n \in \mathbb{Z}^d} \cdots \sum_{z_1 \in \mathbb{Z}^d} \widehat{p}_{t_n,t}(x, z_n) \cdots \widehat{p}_{t_1,t_2}(z_2, z_1) \times$$

$$\times \int_0^{t_1} \sum_{y \in \mathbb{Z}^d} \widehat{p}_{r,t_1}(z_1, y) \, d \left| M_r(\eta_{r^-}[\omega], y)[\omega] \right| \; < \; \infty \, . \qquad (\text{A.18})$$

(3) *For $\mathbb{P}$-a.e. $\omega \in \Omega$, $k = 0, \ldots, n$ and all initial configurations $\eta \in \{0, 1\}^{\mathbb{Z}^d}$, the infinite summation*

$$\sum_{y \in \mathbb{Z}^d} \widehat{p}_{t_k,t_{k+1}}(x, y) \eta_{t_k}(y)[\omega] + \int_{t_k}^{t_{k+1}} \sum_{y \in \mathbb{Z}^d} \widehat{p}_{r,t_{k+1}}(x, y) \, dM_r(\eta_{r^-}[\omega], y)[\omega] \qquad (\text{A.19})$$

*is absolutely convergent and equals $\eta_{t_{k+1}}(x)[\omega]$.*

*Proof.* For item (1), all terms being non-negative, we can reorder the summation on the l.h.s. in (A.16) to obtain

$$\sum_{z_1 \in \mathbb{Z}^d} \cdots \sum_{z_{n-1} \in \mathbb{Z}^d} p_{0,t_1}(x, z_1) \cdots p_{t_{n-2},t_{n-1}}(z_{n-2}, z_{n-1}) \times$$

$$\times \left( \sum_{m_n=0}^{\infty} \sum_{z_n : |z_n - x| = m_n} p_{t_{n-1},t_n}(z_{n-1}, z_n) \sum_{y : |y - x| = m} p_{t_n,t}(z_n, y) \right) . \qquad (\text{A.20})$$

As a consequence of the graphical construction of forward and backward random walks (see Appendix A), triangle inequality and (A.15), we have, for all $z_n \in \mathbb{Z}^d$ such that $|z_n - x| = m_n$, first when $m_n \neq m$,

$$\begin{aligned}
\sum_{y : |y - x| = m} p_{t_n,t}(z_n, y) &= \mathbb{P}\left( |X_{t_n,t}^{z_n} - x| = m \right) \\
&\leq \mathbb{P}\left( \exists y \in \mathbb{Z}^d : |y - x| = m \text{ and } y \in \mathcal{G}_{[t_n,t]}(z_n) \right) \\
&\leq \mathbb{P}\left( \exists w \in \mathbb{Z}^d : |w - z_n| = |m - m_n| \text{ and } w \in \mathcal{G}_{[t_n,t]}(z_n) \right) \\
&\leq e^{-\chi(t-t_n)|m - m_n|} .
\end{aligned} \qquad (\text{A.21})$$

Then when $m_n = m$, we simply bound

$$\sum_{y : |y - x| = m} p_{t_n,t}(z_n, y) \; \leq \; 1 \, . \qquad (\text{A.22})$$

Hence, by (A.21) and (A.22), (A.20) is bounded above by

$$\sum_{m_n=0}^{\infty} e^{-\chi(t-t_n)|m - m_n|} \left( \sum_{z_n : |z_n - x| = m_n} \sum_{z_1 \in \mathbb{Z}^d} \cdots \sum_{z_{n-1} \in \mathbb{Z}^d} p_{0,t_1}(x, z_1) \cdots p_{t_{n-1},t_n}(z_{n-1}, z_n) \right) .$$

By iterating this procedure for a finite number of steps, we obtain the following upper

bound for (A.20):

$$\sum_{m_n=0}^{\infty} \cdots \sum_{m_2=0}^{\infty} e^{-\chi(t-t_n)|m-m_n|} \cdots e^{-\chi(t_2-t_1)|m_2-m_1|} \left( \sum_{m_1=0}^{\infty} \sum_{z_1:|z_1-x|=m_1} p_{0,t_1}(x,z_1) \right).$$

If we bound the last summation in parenthesis as follows (see also (A.21))

$$\sum_{z_1:|z_1-x|=m_1} p_{0,t_1}(x,z_1) \leq \mathbb{P}\left(\exists v \in \mathbb{Z}^d : |v-x| = m_1 \text{ and } v \in \mathcal{G}_{[0,t_1]}(x)\right) \leq e^{-\chi(t_1)m_1},$$

(A.23)

we obtain

$$\sum_{y:|y-x|=m} \sum_{z_1 \in \mathbb{Z}^d} \cdots \sum_{z_n \in \mathbb{Z}^d} p_{0,t_1}(x,z_1) \cdots p_{t_n,t}(z_n,y)$$

$$\leq \sum_{m_n=0}^{\infty} \cdots \sum_{m_2=0}^{\infty} \sum_{m_1=0}^{\infty} e^{-\chi(t-t_n)|m-m_n|} \cdots e^{-\chi(t_2-t_1)|m_2-m_1|} e^{-\chi(t_1)m_1},$$

hence the bound (A.16). An analogous argument yields (A.17).

We now prove item (2). For all initial conditions $\eta \in \{0,1\}^{\mathbb{Z}^d}$, for all realizations $\omega \in \Omega$, by the definition of $dM_r$ in (4.16) and the bound $\eta_{r^-}(z) \leq 1$ for all $z \in \mathbb{Z}^d$, we have that

$$\sum_{z_n \in \mathbb{Z}^d} \cdots \sum_{z_1 \in \mathbb{Z}^d} \widehat{p}_{t_n,t}(x,z_n) \cdots \widehat{p}_{t_1,t_2}(z_2,z_1) \int_0^{t_1} \sum_{y \in \mathbb{Z}^d} \widehat{p}_{r,t_1}(z_1,y) \, d\left|M_r(\eta_{r^-}[\omega],y)[\omega]\right|$$

is bounded above by

$$\sum_{z_n \in \mathbb{Z}^d} \cdots \sum_{z_1 \in \mathbb{Z}^d} \widehat{p}_{t_n,t}(x,z_n) \cdots \widehat{p}_{t_1,t_2}(z_2,z_1) \int_0^{t_1} \sum_{y \in \mathbb{Z}^d} \widehat{p}_{r,t_1}(z_1,y) \sum_{z:\{y,z\} \in E_d} \lambda_r(\{y,z\}) \, dr$$

$$+ \sum_{z_n \in \mathbb{Z}^d} \cdots \sum_{z_1 \in \mathbb{Z}^d} \widehat{p}_{t_n,t}(x,z_n) \cdots \widehat{p}_{t_1,t_2}(z_2,z_1) \int_0^{t_1} \sum_{y \in \mathbb{Z}^d} \widehat{p}_{r,t_1}(z_1,y) \sum_{z:\{y,z\} \in E_d} d\mathcal{N}_r(\{y,z\})[\omega].$$

We estimate the two terms on the above r.h.s. separately. First, by Assumption 2.1 and (A.17), we obtain

$$\sum_{z_n \in \mathbb{Z}^d} \cdots \sum_{z_1 \in \mathbb{Z}^d} \widehat{p}_{t_n,t}(x,z_n) \cdots \widehat{p}_{t_1,t_2}(z_2,z_1) \times$$

$$\times \int_0^{t_1} \sum_{y \in \mathbb{Z}^d} \widehat{p}_{r,t_1}(z_1,y) \sum_{z:\{y,z\} \in E_d} \lambda_r(\{y,z\}) \, dr \leq 2Cd\mathfrak{a}t_1 \sum_{m=0}^{\infty} e^{-\overline{\chi}m} < \infty.$$

For the second term, we observe that, for all $m \in \mathbb{N}$ and $y \in \mathbb{Z}^d$ with $|y - x| = m$, by independence of the Poisson processes over the bonds and Assumption 2.1, we have

$$\mathbb{P}\left(\sum_{z:\{y,z\}\in E_d} \mathcal{N}_{t_1}(\{y,z\}) > m\right) \leq \sum_{k=m+1}^{\infty} \frac{(2d\mathfrak{a}t_1)^k e^{-2d\mathfrak{a}t_1}}{k!} \leq c\,\frac{(2d\mathfrak{a}t_1)^m}{m!}$$

for some constant $c > 0$ independent of $m$. As a consequence, we obtain

$$\sum_{m=0}^{\infty} \sum_{y:|y-x|=m} \mathbb{P}\left(\sum_{z:\{y,z\}\in E_d} \mathcal{N}_{t_1}(\{y,z\}) > m\right)$$
$$\leq \sum_{m=0}^{\infty} \sum_{y:|y-x|=m} c\,\frac{(2d\mathfrak{a}t_1)^m}{m!} \leq \sum_{m=0}^{\infty} c\,\frac{(4d^2\mathfrak{a}t_1)^m}{m!} < \infty.$$

Hence, by a Borel-Cantelli argument, we can conclude that, for $\mathbb{P}$-a.e. $\omega \in \Omega$, there exists a constant $c[\omega] > 0$ for which

$$\sum_{z:\{y,z\}\in E_d} \mathcal{N}_{t_1}(\{y,z\})[\omega] \leq c[\omega]\,m \tag{A.24}$$

holds for all $m \in \mathbb{N}$ and $y \in \mathbb{Z}^d$ with $|y - x| = m$. Therefore, for $\mathbb{P}$-a.e. $\omega \in \Omega$, by (A.24), we get

$$\sum_{z_n \in \mathbb{Z}^d} \cdots \sum_{z_1 \in \mathbb{Z}^d} \widehat{p}_{t_n,t}(x,z_n)\cdots \widehat{p}_{t_1,t_2}(z_2,z_1) \times$$
$$\times \int_0^{t_1} \sum_{y \in \mathbb{Z}^d} \widehat{p}_{r,t_1}(z_1,y) \sum_{z:\{y,z\}\in E_d} \mathrm{d}\mathcal{N}_r(\{y,z\})[\omega]$$
$$\leq \sum_{m=0}^{\infty} \left(\sup_{0\leq r\leq t_1} \sum_{y:|y-x|=m} \sum_{z_n \in \mathbb{Z}^d}\cdots \sum_{z_1 \in \mathbb{Z}^d} \widehat{p}_{t_n,t}(x,z_n)\cdots \widehat{p}_{r,t_1}(z_1,y)\right) 2d\mathfrak{a}t_1 c[\omega]m.$$

The term in parenthesis, by using (A.17), is exponentially small in $m \in \mathbb{N}$, yielding (A.18).

For item (3), for all initial conditions $\eta \in \{0,1\}^{\mathbb{Z}^d}$, we observe that, in view of (A.17), the bound $\eta_0(y) \leq 1$ for all $y \in \mathbb{Z}^d$ and (A.18), for $\mathbb{P}$-a.e. $\omega \in \Omega$, the infinite summation in (A.19) is absolutely convergent. More precisely, for all $\varepsilon > 0$, there exists an integer $n_k = n_{[t_k,t_{k+1}],\varepsilon}(x)[\omega] > 0$ such that

$$\left|\sum_{y:|y-x|\geq n_k} \widehat{p}_{t_k,t_{k+1}}(x,y)\eta_{t_k}(y)[\omega] + \int_{t_k}^{t_{k+1}} \sum_{y:|y-x|\geq n_k} \widehat{p}_{r,t_{k+1}}(x,y)\,\mathrm{d}M_r(\eta_{r^-}[\omega],y)[\omega]\right| < \varepsilon.$$
$$\tag{A.25}$$

Once we have determined $n_k = n_{[t_k, t_{k+1}], \varepsilon}(x)[\omega]$ for which (A.25) is in force, let us define the finite subset $\mathcal{U}_k = \mathcal{U}_{[t_k, t_{k+1}], \varepsilon}(x)[\omega]$ of $\mathbb{Z}^d$ obtained as union of all active islands in $[t_k, t_{k+1}]$ which contain at least a site at a distance $n_k$ (or less) from $x \in \mathbb{Z}^d$. Therefore, for all finite $\mathcal{V} \subset \mathbb{Z}^d$ containing $\mathcal{U}_k$, we have that the absolute value of

$$\left( \sum_{y \in \mathcal{V}} \widehat{p}_{t_k, t_{k+1}}(x, y) \eta_{t_k}(y)[\omega] + \int_{t_k}^{t_{k+1}} \sum_{y \in \mathcal{V}} \widehat{p}_{r, t_{k+1}}(x, y) \, dM_r(\eta_{r^-}[\omega], y)[\omega] \right) - \eta_{t_{k+1}}(x)[\omega]$$

$$= \sum_{y \in \mathcal{V} \setminus \mathcal{U}_k} \widehat{p}_{t_k, t_{k+1}}(x, y) \eta_{t_k}(y)[\omega] + \int_{t_k}^{t_{k+1}} \sum_{y \in \mathcal{V} \setminus \mathcal{U}_k} \widehat{p}_{r, t_{k+1}}(x, y) \, dM_r(\eta_{r^-}[\omega], y)[\omega]$$

(in this identity we used Lemma A.4) is bounded above by

$$\sum_{y \in \mathcal{V} \setminus \mathcal{U}_k} \left( \widehat{p}_{t_k, t_{k+1}}(x, y) \eta_{t_k}(y)[\omega] + \int_{t_k}^{t_{k+1}} \widehat{p}_{r, t_{k+1}}(x, y) \, d|M_r(\eta_{r^-}[\omega], y)[\omega]| \right)$$

$$\leq \sum_{y : |y - x| \geq n_k} \left( \widehat{p}_{t_k, t_{k+1}}(x, y) \eta_{t_k}(y)[\omega] + \int_{t_k}^{t_{k+1}} \widehat{p}_{r, t_{k+1}}(x, y) \, d|M_r(\eta_{r^-}[\omega], y)[\omega]| \right),$$

where this last inequality follows from $\mathcal{V} \setminus \mathcal{U}_k \subset \{y \in \mathbb{Z}^d : |y - x| \geq n_k\}$. By (A.25) the proof is concluded. $\square$

We are now ready to conclude the proof of Proposition 4.3.

*Proof of Proposition 4.3.* For all initial conditions $\eta \in \{0, 1\}^{\mathbb{Z}^d}$, $\mathbb{P}$-a.e. realization $\omega \in \Omega$, $x \in \mathbb{Z}^d$ and $t > 0$, by applying Lemma A.5(3) for all $k = 0, \ldots, n$, and reordering the summations thanks to Lemma A.5(1)–(2),

$$\eta_t(x)[\omega] = \sum_{z_n \in \mathbb{Z}^d} \widehat{p}_{t_n, t}(x, z_n) \eta_{t_n}(z_n)[\omega] + \int_{t_n}^{t} \sum_{z_n \in \mathbb{Z}^d} \widehat{p}_{r, t}(x, z_n) \, dM_r(\eta_{r^-}[\omega], z_n)[\omega]$$

$$= \sum_{z_n \in \mathbb{Z}^d} \cdots \sum_{z_1 \in \mathbb{Z}^d} \sum_{y \in \mathbb{Z}^d} \widehat{p}_{t_n, t}(x, z_n) \cdots \widehat{p}_{0, t_1}(z_1, y) \eta_0(y)$$

$$+ \sum_{z_n \in \mathbb{Z}^d} \cdots \sum_{z_1 \in \mathbb{Z}^d} \widehat{p}_{t_n, t}(x, z_n) \cdots \widehat{p}_{t_1, t_2}(z_2, z_1) \int_0^{t_1} \sum_{y \in \mathbb{Z}^d} \widehat{p}_{r, t_1}(z_1, y) \, dM_r(\eta_{r^-}[\omega], y)[\omega]$$

$$+ \ldots + \int_{t_n}^{t} \sum_{z_n \in \mathbb{Z}^d} \widehat{p}_{r, t}(x, z_n) \, dM_r(\eta_{r^-}[\omega], z_n)[\omega],$$

where the hidden terms are of the following form:

$$\sum_{z_n \in \mathbb{Z}^d} \cdots \sum_{z_k \in \mathbb{Z}^d} \widehat{p}_{t_n, t}(x, z_n) \cdots \widehat{p}_{t_k, t_{k+1}}(z_{k+1}, z_k) \times$$

$$\times \int_{t_{k-1}}^{t_k} \sum_{z_{k-1} \in \mathbb{Z}^d} \widehat{p}_{r,t_k}(z_k, z_{k-1}) \, \mathrm{d}M_r(\eta_{r^-}[\omega], z_{k-1})[\omega] \,, \quad \text{(A.26)}$$

with $k = 1, \ldots, n - 1$. Hence, by using Chapman-Kolmogorov equation (4.6) for the backward transition probabilities, each term in (A.26) equals

$$\int_{t_{k-1}}^{t_k} \sum_{z_{k-1} \in \mathbb{Z}^d} \widehat{p}_{r,t}(x, z_{k-1}) \, \mathrm{d}M_r(\eta_{r^-}[\omega], z_{k-1})[\omega] \,.$$

Thus, by piecing together the above integrals for all $k = 0, \ldots, n$, we finally obtain (4.17).

$\square$

## B  Tightness criterion

We present a tightness criterion for processes in the Skorokhod space $D([0, T], \mathbb{R})$ of real-valued càdlàg functions on $[0, T]$ (see e.g. [5]). This criterion relies on the notion of *uniform conditional stochastic continuity* of a process ([50, Appendix A]). The study of this property allows to extract information on the modulus of continuity of the trajectories.

By following closely the argument in [50, Appendix A], we get a quantitative estimate for the modulus of continuity leading to a sufficient condition for tightness. To the best of our knowledge, this strategy has not been remarked before with this purpose, therefore we provide below a detailed proof.

As a first step, we specify the topological setting as in [5].

**Definition B.1** (MODULUS OF CONTINUITY). *Given $z : [0, T] \to \mathbb{R}$ a bounded function, for all $\delta > 0$, the $\delta$-modulus of continuity $w_z'''(\delta)$ ([5, Problem 12.4, p. 137]) for the function $z$ is given by*

$$w_z'''(\delta) := \max \left\{ \sup_{\substack{0 \leq s \leq t \leq T \\ t-s < \delta}} \inf_{r \in (s,t)} \max \{w_z(s, r), w_z(r, t)\}, \, |z_\delta - z_0|, \, |z_{T^-} - z_{T-\delta}| \right\},$$

*where*

$$w_z(s, t) := \sup_{s \leq s' \leq t' \leq t} |z_{t'} - z_{s'}| \,.$$

Roughly speaking, given $\delta > 0$, the $\delta$-modulus of continuity $w_z'''(\delta)$ (referred to as $\overline{w}_z(\delta)$ in [50, Appendix A.2]) "allows" for one jump in intervals of size at most $\delta$. We refer to [5, Chapter 3. Section 12] for further details on $w'''$ and its relation to the space $D([0, T], \mathbb{R})$. Note that our definition of $w_z'''(\delta)$ slightly differs from the one given in [5, Problem 12.4, p. 137] as we include also information about $z$ near 0 and $T$, i.e. $|z_\delta - z_0|$ and $|z_{T^-} - z_{T-\delta}|$.

In what follows, we state a general tightness criterion, namely Theorem B.2, in $D([0, T], \mathbb{R})$ in terms of the modulus of continuity $w'''$ introduced above. We remark that Theorem B.2 below is a rewriting of Theorem 13.2, the corresponding Corollary and (13.8) to be found at pp. 139–141 of [5]. There the author refers to moduli of continuity ($w'$ and $w''$ defined in (12.6), p. 122, and (12.27), p. 131, respectively) which are different, yet "equivalent" (cf. [5, (12.31)–(12.32), p. 132] and [5, Problem 12.4, p. 137], respectively) to the one we employ.

**Theorem B.2** (TIGHTNESS CRITERION [5, pp. 139–141]). *A family of probability measures $\{P^N, N \in \mathbb{N}\}$ on $D([0, T], \mathbb{R})$, whose canonical coordinate processes are denoted by $\{Z_\cdot^N, N \in \mathbb{N}\}$, is* tight *if the following conditions hold:*

(T1*) *For all $t$ in a dense subset of $[0, T]$ containing $T$,*

$$\lim_{m \to \infty} \limsup_{N \to \infty} P^N \left( |Z_t^N| > m \right) = 0.$$

(T2*) *For all $\varepsilon > 0$,*

$$\lim_{\delta \to 0} \limsup_{N \to \infty} P^N \left( w'''_{Z^N}(\delta) > \varepsilon \right) = 0.$$

In Theorem B.4 below, we will present a condition alternative to (T2*) on the uniform control of the modulus of continuity $w'''$. First we need Theorem B.3 below, which is a slight modification of [50, Theorem A.6]. Indeed, the proof of Theorem B.3 follows closely the one of [50, Theorem A.6]. The major difference lies in the last part: there the two proofs yield different upper bounds (cf. [50, Eq. (A.7)] and (B.2) below). Yet, for the sake of completeness, we include the whole proof at the end of this section.

**Theorem B.3** ([50, THEOREM A.6]). *Let $\{Z_t, t \geq 0\}$ be a continuous-time real-valued stochastic process, whose associated distribution, expectation and filtration are indicated by $P$, $E$ and $\{F_t, t \geq 0\}$, respectively.*

*Fix $T > 0$ and $\varepsilon > 0$ and suppose that there exist a positive value $h_\varepsilon > 0$ and a deterministic function $\psi_\varepsilon : [0, h_\varepsilon] \to [0, \infty)$ such that:*

(i) *$\psi_\varepsilon$ is non-decreasing.*

(ii) *For all $h \in [0, h_\varepsilon]$ and $t \in [0, T]$, we have*

$$P\left( |Z_{t+h} - Z_t| > \varepsilon \mid F_t \right) \leq \psi_\varepsilon(h), \quad \text{P-a.s.} \,. \qquad (B.1)$$

*Then the following bound on the modulus of continuity $w'''_Z$*

$$P\left( w'''_Z(h) > 4\varepsilon \right) \leq 2(k+1)\psi_\varepsilon(h) + 4\psi_\varepsilon\left(\tfrac{2T}{k}\right) \qquad (B.2)$$

*holds for all $h \in [0, h_\varepsilon]$ and $k > k_\varepsilon := 2T/h_\varepsilon$.*

We remark that, if both (B.2) and $\psi_\varepsilon(h) \to 0$ as $h \to 0$ hold for all $\varepsilon > 0$, then the process $\{Z_t, t \in [0, T]\}$ can be realized in the Skorokhod space $D([0, T], \mathbb{R})$ (see [50, Theorem A.6] for further details).

**Theorem B.4** (UNIFORM CONDITIONAL STOCHASTIC EQUICONTINUITY & (T2*))**.** *Let* $\{P^N, N \in \mathbb{N}\}$ *and* $\{Z^N_\cdot, N \in \mathbb{N}\}$ *be as in Theorem B.2. Let, for all* $N \in \mathbb{N}$, $\{Z^N_t, 0 \leq t \leq T\}$ *be adapted to the filtration* $\{F^N_t, 0 \leq t \leq T\}$.

*Fix* $\varepsilon > 0$ *and suppose that there exist values* $h_\varepsilon > 0$ *and* $N_\varepsilon \in \mathbb{N}$ *such that for all* $N \geq N_\varepsilon$ *there exist deterministic functions* $\psi^N_\varepsilon, \psi_\varepsilon : [0, h_\varepsilon] \to [0, \infty)$ *and non-negative values* $\phi^N_\varepsilon$ *satisfying, for all* $N \geq N_\varepsilon$, *the following properties:*

(i) *The functions* $\psi^N_\varepsilon$ *are non-decreasing.*

(ii) *For all* $h \in [0, h_\varepsilon]$ *and* $t \in [0, T]$, *we have*
$$P^N \left( |Z^N_{t+h} - Z^N_t| > \varepsilon \mid F^N_t \right) \leq \psi^N_\varepsilon(h), \quad P^N\text{-a.s. .}$$

(iii) *For all* $h \in [0, h_\varepsilon]$, *we have* $\psi^N_\varepsilon(h) \leq \psi_\varepsilon(h) + \phi^N_\varepsilon$.

(iv) $\phi^N_\varepsilon \to 0$ *as* $N \to \infty$.

(v) $\psi_\varepsilon(h) \to 0$ *as* $h \to 0$.

*Then we obtain*
$$\lim_{\delta \to 0} \limsup_{N \to \infty} P^N \left( w'''_{Z^N}(\delta) > 4\varepsilon \right) = 0. \tag{B.3}$$

*If this is true for all* $\varepsilon > 0$, *then condition* (T2*) *in Theorem B.2 holds for* $\{Z^N_\cdot, N \in \mathbb{N}\}$.

*Proof.* Fix $\varepsilon > 0$. Due to (i) and (ii) we can apply Theorem B.3 to get an estimate for $P^N(w'''_{Z^N}(h) > 4\varepsilon)$ of the form (B.2) with $\psi^N_\varepsilon$. By using, in addition, (iii), we obtain the bound
$$P^N \left( w'''_{Z^N}(h) > 4\varepsilon \right) \leq 2(k+1)\left(\psi_\varepsilon(h) + \phi^N_\varepsilon\right) + 4\left(\psi_\varepsilon(\tfrac{2T}{k}) + \phi^N_\varepsilon\right),$$
which is valid for all $h \in [0, h_\varepsilon]$, $N \geq N_\varepsilon$ and $k > 2T/h_\varepsilon$. Now observe that, by (iv), we have
$$\limsup_{N \to \infty} P^N \left( w'''_{Z^N}(h) > 4\varepsilon \right) \leq 2(k+1)\psi_\varepsilon(h) + 4\psi_\varepsilon(\tfrac{2T}{k}). \tag{B.4}$$

We are left to show that the r.h.s. in (B.4) vanishes as $h \to 0$. We use the fact that $\psi_\varepsilon(h) \to 0$ as $h \to 0$ in (v). Indeed, by first taking the limit as $h \to 0$, we get
$$\lim_{h \to 0} \limsup_{N \to \infty} P^N \left( w'''_{Z^N}(h) > 4\varepsilon \right) \leq 4\psi_\varepsilon(\tfrac{2T}{k}),$$
which holds for all $k > k_\varepsilon = 2T/h_\varepsilon$. A further limit with $k \to \infty$ on the r.h.s. above and (v) yield (B.3). □

*Proof of Theorem B.3.* We follow here [50, Theorem A.6] and refer, without further mention, to the process stopped at time $T$ when writing $\{Z_t, t \geq 0\}$ in this proof. We fix $\varepsilon > 0$, $\tau_{\varepsilon,0} = 0$ and define $\tau_{\varepsilon,1}$ as the first time $|Z_t - Z_0|$ exceeds $2\varepsilon$, $\tau_{\varepsilon,1} + \tau_{\varepsilon,2}$ as the first time $t > \tau_{\varepsilon,1}$ for which $|Z_t - Z_{\tau_{\varepsilon,1}}|$ does and so on, up to reach time $T$ and with the convention that, if $\tau_{\varepsilon,1} + \ldots + \tau_{\varepsilon,n} > T$, we set $\tau_{\varepsilon,1} + \ldots + \tau_{\varepsilon,n}$ equal to $T + 2\, h_\varepsilon\, n$. As a consequence of these definitions, if we define $\sigma_{\varepsilon,n} = \tau_{\varepsilon,1} + \ldots + \tau_{\varepsilon,n}$, we have: for all $n \in \mathbb{N}_0$,

$$P\left( \sup_{\sigma_{\varepsilon,n} \leq s \leq t < \sigma_{\varepsilon,n+1}} |Z_t - Z_s| \leq 4\varepsilon \right) = 1\, ,$$

and, for all $h \in [0, h_\varepsilon]$,

$$P\left( \sup_{0 \leq h' \leq h} |Z_{\sigma_{\varepsilon,n}+h'} - Z_{\sigma_{\varepsilon,n}}| > 2\varepsilon \,\Big|\, F_{\sigma_{\varepsilon,n}} \right) = P\left( \tau_{\varepsilon,n+1} \leq h \,|\, F_{\sigma_{\varepsilon,n}} \right)\, . \quad \text{(B.5)}$$

We rewrite the probability in (B.5) as follows:

$$\begin{aligned} P\left( \tau_{\varepsilon,n+1} \leq h \,|\, F_{\sigma_{\varepsilon,n}} \right) &= P\left( \tau_{\varepsilon,n+1} \leq h,\, |Z_{\sigma_{\varepsilon,n}+h} - Z_{\sigma_{\varepsilon,n}}| \leq \varepsilon \,|\, F_{\sigma_{\varepsilon,n}} \right) \\ &\quad + P\left( \tau_{\varepsilon,n+1} \leq h,\, |Z_{\sigma_{\varepsilon,n}+h} - Z_{\sigma_{\varepsilon,n}}| > \varepsilon \,|\, F_{\sigma_{\varepsilon,n}} \right)\, . \end{aligned} \quad \text{(B.6)}$$

Concerning the first term on the r.h.s. in (B.6), we have the following upper bound (recall that $\sigma_{\varepsilon,n+1} = \sigma_{\varepsilon,n} + \tau_{\varepsilon,n+1}$):

$$P\left( \tau_{\varepsilon,n+1} \leq h,\, |Z_{\sigma_{\varepsilon,n}+h} - Z_{\sigma_{\varepsilon,n+1}}| > \varepsilon \,|\, F_{\sigma_{\varepsilon,n}} \right)\, ,$$

which, in turn, rewrites as follows:

$$\begin{aligned} &E\left[ E\left[ \mathbf{1}_{\{\tau_{\varepsilon,n+1} \leq h\}} \mathbf{1}_{\{|Z_{\sigma_{\varepsilon,n}+h} - Z_{\sigma_{\varepsilon,n+1}}| > \varepsilon\}} \,\Big|\, F_{\sigma_{\varepsilon,n+1}} \right] \,\Big|\, F_{\sigma_{\varepsilon,n}} \right] \\ &= E\left[ \mathbf{1}_{\{\tau_{\varepsilon,n+1} \leq h\}}\, E\left[ \mathbf{1}_{\{|Z_{\sigma_{\varepsilon,n}+h} - Z_{\sigma_{\varepsilon,n+1}}| > \varepsilon\}} \,\Big|\, F_{\sigma_{\varepsilon,n+1}} \right] \,\Big|\, F_{\sigma_{\varepsilon,n}} \right]\, . \end{aligned}$$

By (B.1) - which holds true also when considering $\sigma$-fields associated to stopping times being the bound (B.1) uniform in time - we obtain, P-a.s.,

$$\begin{aligned} &P\left( \tau_{\varepsilon,n+1} \leq h,\, |Z_{\sigma_{\varepsilon,n}+h} - Z_{\sigma_{\varepsilon,n}}| \leq \varepsilon \,|\, F_{\sigma_{\varepsilon,n}} \right) \\ &\leq P\left( \tau_{\varepsilon,n+1} \leq h,\, |Z_{\sigma_{\varepsilon,n}+h} - Z_{\sigma_{\varepsilon,n}+\tau_{\varepsilon,n+1}}| > \varepsilon \,|\, F_{\sigma_{\varepsilon,n}} \right) \\ &\leq E\left[ \mathbf{1}_{\{\tau_{\varepsilon,n+1} \leq h\}}\, \psi_\varepsilon(h - \tau_{\varepsilon,n+1}) \,|\, F_{\sigma_{\varepsilon,n}} \right] \\ &\leq \psi_\varepsilon(h)\, P\left( \tau_{\varepsilon,n+1} \leq h \,|\, F_{\sigma_{\varepsilon,n}} \right)\, , \end{aligned} \quad \text{(B.7)}$$

where in the last inequality we used the monotonicity of $\psi_\varepsilon$ ($\psi_\varepsilon(h') \leq \psi_\varepsilon(h'')$ if $h' \leq h''$). Analogously, being the bound (B.1) uniform in time, for the second term on the r.h.s. in (B.6), we have, for all $h \in [0, h_\varepsilon]$ and P-a.s. as a consequence of (B.1), the bound in (B.1) being uniform in time,

$$P\left( \tau_{\varepsilon,n+1} \leq h,\, |Z_{\sigma_{\varepsilon,n}+h} - Z_{\sigma_{\varepsilon,n}}| > \varepsilon \,|\, F_{\sigma_{\varepsilon,n}} \right)$$

$$\leq \; \mathrm{P}\left(|Z_{\sigma_{\varepsilon,n}+h} - Z_{\sigma_{\varepsilon,n}}| > \varepsilon \,|\, \mathrm{F}_{\sigma_{\varepsilon,n}}\right) \;\leq\; \psi_\varepsilon(h) \,. \tag{B.8}$$

As a consequence of (B.5), (B.6), (B.7) and (B.8), we obtain, for all $h \in [0, h_\varepsilon]$ and $n \in \mathbb{N}_0$,

$$\mathrm{P}\left(\sup_{0 \leq h' \leq h} |Z_{\sigma_{\varepsilon,n}+h'} - Z_{\sigma_{\varepsilon,n}}| > 2\varepsilon \,\bigg|\, \mathrm{F}_{\sigma_{\varepsilon,n}}\right) \;\leq\; 2\psi_\varepsilon(h), \quad \text{P-a.s.} \,. \tag{B.9}$$

Recall Definition B.1 of the modulus of continuity $w'''$. For any choice of $k \in \mathbb{N}$, the probability $\mathrm{P}\left(w_Z'''(h) > 4\varepsilon\right)$ can be bounded above by

$$\mathrm{P}\left(w_Z'''(h) > 4\varepsilon,\, \sigma_{\varepsilon,k} > T,\, \min\{\tau_{\varepsilon,1}, \ldots, \tau_{\varepsilon,k}\} > h,\, |Z_{T^-} - Z_{T-h}| \leq 4\varepsilon\right) \tag{B.10}$$

$$+ \; \mathrm{P}\left(\sigma_{\varepsilon,k} \leq T\right) \;+\; \mathrm{P}\left(\min\{\tau_{\varepsilon,1}, \ldots, \tau_{\varepsilon,k}\} \leq h\right) \;+\; \mathrm{P}\left(|Z_{T^-} - Z_{T-h}| > 4\varepsilon\right) \,. \tag{B.11}$$

The probability in (B.10) vanishes. Indeed, if the events

$$\{\sigma_{\varepsilon,k} > T\} \;\text{ and }\; \{\min\{\tau_{\varepsilon,1}, \ldots, \tau_{\varepsilon,k}\} > h\}$$

occur, then necessarily in any subinterval of size $h$ of $[0, T]$ there can be at most one $\sigma_{\varepsilon,\ell}$, for some $0 \leq \ell \leq k$, making, together with

$$\{|Z_{T^-} - Z_{T-h}| \leq 4\varepsilon\} \,,$$

the event $\{w_Z'''(h) > 4\varepsilon\}$ impossible.

Now we estimate each term in (B.11) and consider $h \in [0, h_\varepsilon]$. For the second one, by (B.5) and (B.9), we get

$$\mathrm{P}\left(\min\{\tau_{\varepsilon,1}, \ldots, \tau_{\varepsilon,k}\} \leq h\right) \;\leq\; \sum_{\ell=1}^{k} \mathrm{E}\left[\mathrm{P}\left(\tau_{\varepsilon,\ell} \leq h \,|\, \mathrm{F}_{\sigma_{\varepsilon,\ell-1}}\right)\right] \;\leq\; 2k\psi_\varepsilon(h) \,. \tag{B.12}$$

For the third term, we obtain

$$\mathrm{P}\left(|Z_{T^-} - Z_{T-h}| > 4\varepsilon\right) \;\leq\; \mathrm{P}\left(\sup_{T-h \leq t < T} |Z_t - Z_{T-h}| > 2\varepsilon\right)$$

$$= \; \mathrm{E}\left[\mathrm{P}\left(\sup_{T-h \leq t < T} |Z_t - Z_{T-h}| > 2\varepsilon \,\bigg|\, \mathrm{F}_{T-h}\right)\right] \;\leq\; 2\psi_\varepsilon(h), \tag{B.13}$$

where in the last inequality we argued as to obtain (B.9) and used (B.1). It is slightly more involved to control the first term in (B.11). We have, for all $\delta \in [0, h_\varepsilon]$,

$$T\,\mathrm{P}\left(\sigma_{\varepsilon,k} \leq T\right) \;\geq\; \mathrm{E}\left[\sigma_{\varepsilon,k}\mathbf{1}_{\{\sigma_{\varepsilon,k} \leq T\}}\right] \;=\; \sum_{\ell=1}^{k} \mathrm{E}\left[\tau_{\varepsilon,\ell}\mathbf{1}_{\{\sigma_{\varepsilon,k} \leq T\}}\right]$$

$$\geq \; \sum_{\ell=1}^{k} \mathrm{E}\left[\tau_{\varepsilon,\ell}\mathbf{1}_{\{\sigma_{\varepsilon,k} \leq T\}}\mathbf{1}_{\{\tau_{\varepsilon,\ell} > \delta\}}\right] \;\geq\; \delta \sum_{\ell=1}^{k} \mathrm{P}\left(\sigma_{\varepsilon,k} \leq T,\, \tau_{\varepsilon,\ell} > \delta\right)$$

$$\geq \delta \sum_{\ell=1}^{k} \mathrm{P}\left(\sigma_{\varepsilon,k} \leq T\right) - \delta \sum_{\ell=1}^{k} \mathrm{P}\left(\tau_{\varepsilon,\ell} \leq \delta\right)$$
$$\geq \delta\, k\, \mathrm{P}\left(\sigma_{\varepsilon,k} \leq T\right) - \delta\, k\, 2\,\psi_{\varepsilon}(\delta)\,.$$

where this last inequality follows from (B.5) and (B.9) as in (B.12). Hence, whenever $\delta k > T$, we obtain

$$\mathrm{P}\left(\sigma_{\varepsilon,k} \leq T\right) \;\leq\; \frac{\delta\, k}{\delta\, k - T}\, 2\,\psi_{\varepsilon}(\delta)\,. \tag{B.14}$$

To conclude, the bounds (B.12), (B.14) with the choice $\delta = \frac{2T}{k}$ (and, as a consequence, $k > 2T/h_{\varepsilon}$) and (B.13) lead to the final result (B.2).

□

## C From the invariance principle to the arbitrary starting point invariance principle in dynamic environment: the uniformly elliptic case

In this appendix, which continues the discussion started in Section 3.1, we show that the invariance principle for the random walk starting at the origin (condition (b$_1$) in Theorem 3.4) combined with the additional assumption of uniform ellipticity of the environment (condition (b$_2$) in Theorem 3.4, a stronger condition than the one contained in 2.1) yields the invariance principle with arbitrary starting points (condition (b) in Theorem 3.2). In particular, this shows that all quenched results in e.g. [1], [3], i.e. (3.16) below, imply (b) in Theorem 3.2, provided that the environment satisfies not only a uniform upper bound (Assumption 2.1), but also a uniform lower bound. This result is the content of Theorem 3.4.

In what follows, we present the proof of Theorem 3.4, whose main steps follow the proof in Appendix A.2 of [10], but here adjusted to the dynamic context. The idea consists in letting again random walks and Brownian motions starting from the origin at time zero hit $\varepsilon$-ball around macroscopic points within a certain positive time ($T > 0$ in [10], $\vartheta > 0$ here). The main difference here w.r.t. [10] is that we send not only $\varepsilon$ – the radius of the ball – to zero, but also $\vartheta$ – roughly speaking, the last time available for the first visit of this $\varepsilon$-ball. Key ingredient is the space-time Hölder continuity of the random walk semigroups obtained from bounds for uniformly elliptic time-dependent conductances taken from [26, Appendix B]. Similar estimates are employed in [10], though obtained there via subsequent approximations of the actual random walk by random walks killed on exiting growing macroscopic balls, cf. [10, p. 1639].

The two main technical steps of the proof of Theorem 3.4 are enclosed in two lemmas, one – Lemma C.3 – showing tightness of the sequence of rescaled random walks and the other one – Lemma C.2 – proving convergence of finite-dimensional distributions. The proofs of both lemmas rely on estimates of transition probabilities and semigroups corresponding to time-inhomogeneous random walks in presence of uniformly elliptic conductances. These estimates have been obtained in various forms by several authors in both static and dynamic contexts, see e.g. [48] and [26], [12], respectively. The estimates we employ are those taken from [26, Appendix B] and, as we will repeatedly refer to them, we display them below for the convenience of the reader.

**Proposition C.1** ([26, APPENDIX B]). *Under condition* (b$_2$) *of uniform ellipticity of the conductances* $\lambda$, *the corresponding random walk transition probabilities and semigroups satisfy the following bounds:*

(I) (UPPER BOUND ON THE KERNEL ([26, PROPOSITION B.3])). *There exists a constant* $C \in (1, \infty)$, *depending only on* $d, \mathfrak{a}$ *and* $\mathfrak{b}$, *such that*

$$p^N_{sN^2, tN^2}(x, y) \leq \frac{C}{N^d (\frac{1}{N} \vee \sqrt{t-s})^d} e^{-\frac{|\frac{x}{N} - \frac{y}{N}|}{\frac{1}{N} \vee \sqrt{t-s}}}, \qquad \text{(C.1)}$$

*for all* $x, y \in \mathbb{Z}^d$, $t \geq s$ *and* $N \in \mathbb{N}$.

(II) (NASH CONTINUITY ESTIMATE ([26, PROPOSITION B.6])). *There exists* $\gamma > 0$ *and* $c > 0$, *depending only on* $d, \mathfrak{a}$ *and* $\mathfrak{b}$, *such that, for every* $N \in \mathbb{N}$ *and* $f \in \ell^\infty(\frac{\mathbb{Z}^d}{N})$,

$$\left| S^N_{sN^2, (t+h)N^2} f(\tfrac{x}{N}) - S^N_{sN^2, tN^2} f(\tfrac{y}{N}) \right| \leq c \sup_{x \in \mathbb{Z}^d} |f(\tfrac{x}{N})| \left( \frac{\sqrt{h} \vee |\tfrac{x}{N} - \tfrac{y}{N}|}{\sqrt{t-s}} \right)^\gamma, \quad \text{(C.2)}$$

*for all* $x, y \in \mathbb{Z}^d$, $h \geq 0$ *and* $t \geq s$.

Given tightness – whose proof is postponed to Lemma C.3 below – of the sequence of rescaled random walks

$$\left\{ \frac{X^{x_N}_{0, tN^2}}{N} : t \in [0, T] \right\}$$

with arbitrary starting points, we complete the proof of Theorem 3.4 by showing convergence of finite-dimensional distributions. The general idea of the proof comes from [10, Appendix A.2], where an arbitrary starting point quenched invariance principle was proved in presence of static conductances over bounded supercritical percolation clusters in any dimension $d \geq 2$ (see also [23] for a second application of this method in presence of static site inhomogeneities).

In what follows we adapt this strategy to the dynamic context, postponing the main technical part – a suitable Hölder continuity, both in space and time, of random walk semigroups – to Lemma C.2 below.

*Proof of Theorem 3.4.* We fix $T > 0$, $u \in \mathbb{R}^d$ and $\{x_N : N \in \mathbb{N}\} \subset \mathbb{Z}^d$ as in the statement of the theorem. While the sequence of rescaled random walks (C.15) is tight by Lemma C.3 below, to establish convergence of finite-dimensional distributions we follow [10, Appendix A.2], as below:

*Step 1.* For any $\varepsilon > 0$, we define the following hitting times of Euclidean $\varepsilon$-balls around $u \in \mathbb{R}^d$ (*Notation:* "$\mathcal{B}(u, \varepsilon)$") of both random walks and Brownian motion started

from the origin:

$$\tau_\varepsilon^N := \inf\left\{t \geq 0 : \frac{X_{0,t}^0}{N} \in \mathcal{B}(u,\varepsilon) \subset \mathbb{R}^d\right\}$$

$$\tau_\varepsilon^\Sigma := \inf\left\{t \geq 0 : B_t^\Sigma \in \mathcal{B}(u,\varepsilon) \subset \mathbb{R}^d\right\}.$$

As a consequence of the invariance principle at the origin (3.16) and continuity of Brownian motion's trajectories, we get the following convergence in law:

$$\frac{\tau_\varepsilon^N}{N^2} \underset{N\to\infty}{\Longrightarrow} \tau_\varepsilon^\Sigma. \tag{c.3}$$

*Step 2.* We fix two deterministic constants $\Theta > \vartheta > 0$. As a consequence of (3.16) on the time interval $[0, T + \Theta]$ and (c.3), we have that the invariance principle holds also for the processes "conditioned to hit $\mathcal{B}(u,\varepsilon)$ before time $\vartheta$ and observed after the hitting"(see also [10] for an analogous statement), namely

$$\left\{\frac{X_{0,\tau_\varepsilon^N+tN^2}^0}{N} \,\bigg|\, \frac{\tau_\varepsilon^N}{N^2} < \vartheta : t \in [0, T]\right\} \underset{N\to\infty}{\Longrightarrow} \left\{B_{\tau_\varepsilon^\Sigma+t}^\Sigma \,|\, \tau_\varepsilon^\Sigma < \vartheta : t \in [0, T]\right\}.$$

As a consequence, for all $n \in \mathbb{N}, 0 \leq t_1 < \ldots < t_n \leq T$ and functions $f_1,\ldots,f_n \in C_0(\mathbb{R}^d)$, we have

$$\mathbb{E}\left[f_1\left(\frac{X_{0,\tau_\varepsilon^N+t_1 N^2}^0}{N}\right)\cdots f_n\left(\frac{X_{0,\tau_\varepsilon^N+t_n N^2}^0}{N}\right) \,\bigg|\, \frac{\tau_\varepsilon^N}{N^2} < \vartheta\right]$$

$$\underset{N\to\infty}{\longrightarrow} \mathsf{E}\left[f_1\left(B_{\tau_\varepsilon^\Sigma+t_1}^\Sigma\right)\cdots f_n\left(B_{\tau_\varepsilon^\Sigma+t_n}^\Sigma\right) \,|\, \tau_\varepsilon^\Sigma < \vartheta\right]. \tag{c.4}$$

*Step 3.* As a consequence of the (strong) Markov property of both processes and (c.4), we obtain

$$\int_0^\vartheta \sum_{\frac{y_N}{N} \in \mathcal{B}(u,\varepsilon)} \mathbb{E}\left[f_1\left(\frac{X_{sN^2,(s+t_1)N^2}^{y_N}}{N}\right)\cdots f_n\left(\frac{X_{sN^2,(s+t_n)N^2}^{y_N}}{N}\right)\right] \mathbb{Q}_{\varepsilon,\vartheta}^N(y_N, \mathrm{d}s)$$

$$\underset{N\to\infty}{\longrightarrow} \int_0^\vartheta \int_{\overline{\mathcal{B}(u,\varepsilon)}} \mathsf{E}\left[f_1\left(B_{(s+t_1)-s}^\Sigma + v\right)\cdots f_n\left(B_{(s+t_n)-s}^\Sigma + v\right)\right] \mathbb{Q}_{\varepsilon,\theta}^\Sigma(\mathrm{d}v, \mathrm{d}s)$$

$$= \int_{\overline{\mathcal{B}(u,\varepsilon)}} \mathsf{E}\left[f_1\left(B_{t_1}^\Sigma + v\right)\cdots f_n\left(B_{t_n}^\Sigma + v\right)\right] \mathbb{Q}_{\varepsilon,\theta}^\Sigma(\mathrm{d}v), \tag{c.5}$$

where $\overline{\mathcal{B}(u,\varepsilon)}$ is the closure of $\mathcal{B}(u,\varepsilon)$,

$$\mathbb{Q}_{\varepsilon,\vartheta}^{N}(y_N, ds) := \mathbb{P}\left(\frac{X_{0,sN^2}^{0}}{N} = \frac{y_N}{N},\ \tau_\varepsilon^N \in ds\ \bigg|\ \frac{\tau_\varepsilon^N}{N^2} < \vartheta\right)$$

is the probability that the random walk $X_{0,\cdot}^0$ started from the origin at time zero has entered the (macroscopic) ball $N \cdot \mathcal{B}(u,\varepsilon)$ at position $y_N \in \mathbb{Z}^d$ at times $r \in (sN^2, (s+\delta)N^2)$ for $\delta \ll 1$, conditioned on $\tau_\varepsilon^N/N^2 < \vartheta$ (and hence, necessarily, $s < \vartheta$) and, analogously, for the Brownian motion:

$$\mathbb{Q}_{\varepsilon,\theta}^{\Sigma}(dv, ds) := \mathsf{P}\left(B_{\tau_\varepsilon^\Sigma}^{\Sigma} \in dv,\ \tau_\varepsilon^\Sigma \in ds\ \big|\ \tau_\varepsilon^\Sigma < \vartheta\right)$$

and its marginal probability

$$\mathbb{Q}_{\varepsilon,\theta}^{\Sigma}(dv) := \int_0^{\theta} \mathbb{Q}_{\varepsilon,\theta}^{\Sigma}(dv, ds).$$

We note that the last identity in (C.5) is due to the time-homogeneity of the law of the Brownian motion.

*Step 4.* By the uniform continuity of the functions $f_1, \ldots, f_n \in C_0(\mathbb{R}^d)$, we have

$$\lim_{\varepsilon \to 0} \limsup_{\vartheta \to 0} \left| \int_{\mathcal{B}(u,\varepsilon)} \mathsf{E}\left[f_1\left(B_{t_1}^\Sigma + v\right) \cdots f_n\left(B_{t_n}^\Sigma + v\right)\right] \mathbb{Q}_{\varepsilon,\theta}^\Sigma(dv) \right. $$
$$\left. - \mathsf{E}\left[f_1\left(B_{t_1}^\Sigma + u\right) \cdots f_n\left(B_{t_n}^\Sigma + u\right)\right] \right| = 0. \quad (C.6)$$

Indeed,

$$\left| \int_{\mathcal{B}(u,\varepsilon)} \mathsf{E}\left[f_1\left(B_{t_1}^\Sigma + v\right) \cdots f_n\left(B_{t_n}^\Sigma + v\right)\right] \mathbb{Q}_{\varepsilon,\theta}^\Sigma(dv) - \mathsf{E}\left[f_1\left(B_{t_1}^\Sigma + u\right) \cdots f_n\left(B_{t_n}^\Sigma + u\right)\right] \right|$$
$$\leq \sup_{\substack{v \in \mathbb{R}^d \\ |u-v| < 2\varepsilon}} \left| \mathsf{E}\left[f_1\left(B_{t_1}^\Sigma + v\right) \cdots f_n\left(B_{t_n}^\Sigma + v\right) - f_1\left(B_{t_1}^\Sigma + u\right) \cdots f_n\left(B_{t_n}^\Sigma + u\right)\right] \right|$$
$$\leq \|f_1\|_\infty \cdots \|f_n\|_\infty \sum_{i=1}^n \sup_{\substack{u,v \in \mathbb{R}^d \\ |u-v| < 2\varepsilon}} |f_i(u) - f_i(v)|.$$

The independence of $\theta$ of the last r.h.s. in the expression above and the aforementioned uniform continuity of the (bounded) functions $f_1, \ldots, f_n \in C_0(\mathbb{R}^d)$ yield (C.6).

*Step 5.* Because of the continuity and the boundedness of the functions $f_1, \ldots, f_n \in C_0(\mathbb{R}^d)$, we apply Lemma C.2 below and obtain

$$\lim_{\varepsilon \to 0} \limsup_{\vartheta \to 0} \limsup_{N \to \infty} \left| \mathbb{E}\left[f_1\left(\frac{X_{0,t_1 N^2}^{x_N}}{N}\right) \cdots f_n\left(\frac{X_{0,t_n N^2}^{x_N}}{N}\right)\right] \right.$$

$$-\int_0^\vartheta \sum_{\frac{y_N}{N}\in\mathcal{B}(u,\varepsilon)} \mathbb{E}\left[f_1\left(\frac{X^{y_N}_{sN^2,(s+t_1)N^2}}{N}\right)\cdots f_n\left(\frac{X^{y_N}_{sN^2,(s+t_n)N^2}}{N}\right)\right]\mathbb{Q}^N_{\varepsilon,\vartheta}(y_N,\mathrm{d}s)\Bigg| = 0.$$
(C.7)

In particular, we note that no property of $\mathbb{Q}^N_{\vartheta,\varepsilon}(\cdot,\cdot)$ other than the fact that it is a probability measure plays any role in the proof of (C.7) above.

*Step 6.* Then we are done, because we get that

$$\lim_{N\to\infty}\left|\mathbb{E}\left[f_1\left(\frac{X^{x_N}_{0,t_1N^2}}{N}\right)\cdots f_n\left(\frac{X^{x_N}_{0,t_nN^2}}{N}\right)\right] - \mathsf{E}\left[f_1\left(B^\Sigma_{t_1}+u\right)\cdots f_n\left(B^\Sigma_{t_n}+u\right)\right]\right|$$

is bounded above by the sum of the following three terms:

$$\lim_{\varepsilon\to 0}\limsup_{\vartheta\to 0}\limsup_{N\to\infty}\left|\mathbb{E}\left[f_1\left(\frac{X^{x_N}_{0,t_1N^2}}{N}\right)\cdots f_n\left(\frac{X^{x_N}_{0,t_nN^2}}{N}\right)\right]\right.$$
$$\left.-\int_0^\vartheta\sum_{\frac{y_N}{N}\in\mathcal{B}(u,\varepsilon)}\mathbb{E}\left[f_1\left(\frac{X^{y_N}_{sN^2,(s+t_1)N^2}}{N}\right)\cdots f_n\left(\frac{X^{y_N}_{sN^2,(s+t_n)N^2}}{N}\right)\right]\mathbb{Q}^N_{\varepsilon,\vartheta}(y_N,\mathrm{d}s)\right|,$$

$$\lim_{N\to\infty}\left|\int_0^\vartheta\sum_{\frac{y_N}{N}\in\mathcal{B}(u,\varepsilon)}\mathbb{E}\left[f_1\left(\frac{X^{y_N}_{sN^2,(s+t_1)N^2}}{N}\right)\cdots f_n\left(\frac{X^{y_N}_{sN^2,(s+t_n)N^2}}{N}\right)\right]\mathbb{Q}^N_{\varepsilon,\vartheta}(y_N,\mathrm{d}s)\right.$$
$$\left.-\int_{\mathcal{B}(u,\varepsilon)}\mathsf{E}\left[f_1\left(B^\Sigma_{t_1}+v\right)\cdots f_n\left(B^\Sigma_{t_n}+v\right)\right]\mathbb{Q}^\Sigma_{\varepsilon,\theta}(\mathrm{d}v)\right|$$

and

$$\lim_{\varepsilon\to 0}\limsup_{\vartheta\to 0}\left|\int_{\mathcal{B}(u,\varepsilon)}\mathsf{E}\left[f_1\left(B^\Sigma_{t_1}+v\right)\cdots f_n\left(B^\Sigma_{t_n}+v\right)\right]\mathbb{Q}^\Sigma_{\varepsilon,\theta}(\mathrm{d}v)\right.$$
$$\left.-\mathsf{E}\left[f_1\left(B^\Sigma_{t_1}+u\right)\cdots f_n\left(B^\Sigma_{t_n}+u\right)\right]\right|.$$

These three last terms, by (C.7), (C.5) and (C.6) respectively, equal zero.

□

In the following proposition we fill the gap in the proof of Theorem 3.4 (*Step 5*).

**Lemma C.2** (HÖLDER CONTINUITY OF RW-SEMIGROUPS). *For all $n \in \mathbb{N}$, $0 < t_1 < \ldots < t_n \leq T$ and continuous and bounded functions $f_1,\ldots,f_n : \mathbb{R}^d \to \mathbb{R}$, we have*

$$\lim_{\varepsilon\to 0}\limsup_{\vartheta\to 0}\limsup_{N\to\infty}\sup_{0\leq h<\vartheta}\sup_{|\frac{x}{N}-\frac{y}{N}|<\varepsilon}$$
(C.8)

$$\left| \mathbb{E}\left[ f_1\left(\frac{X^x_{hN^2,(t_1+h)N^2}}{N}\right) \cdots f_n\left(\frac{X^x_{hN^2,(t_n+h)N^2}}{N}\right) \right] - \mathbb{E}\left[ f_1\left(\frac{X^y_{0,t_1N^2}}{N}\right) \cdots f_n\left(\frac{X^y_{0,t_nN^2}}{N}\right) \right] \right| = 0 \,.$$

*Proof.* We first consider the case $n = 1$:

$$\left| \mathbb{E}\left[ f\left(\frac{X^x_{hN^2,(t+h)N^2}}{N}\right) \right] - \mathbb{E}\left[ f\left(\frac{X^y_{0,tN^2}}{N}\right) \right] \right| = \left| S^N_{hN^2,(t+h)N^2} f(\tfrac{x}{N}) - S^N_{0,tN^2} f(\tfrac{y}{N}) \right| ,$$

which, by the Markov property (Proposition A.2(6)), we bound above by

$$\left| S^N_{hN^2,(t+h)N^2} f(\tfrac{x}{N}) - S^N_{hN^2,tN^2} f(\tfrac{y}{N}) \right| + \left| S^N_{hN^2,tN^2} f(\tfrac{y}{N}) - S^N_{0,hN^2} S^N_{hN^2,tN^2} f(\tfrac{y}{N}) \right| . \quad \text{(c.9)}$$

For the first term in (c.9), we apply Nash continuity estimate (c.2) to obtain

$$\left| S^N_{hN^2,(t+h)N^2} f(\tfrac{x}{N}) - S^N_{hN^2,tN^2} f(\tfrac{y}{N}) \right| \leq c \left( \sup_{z \in \mathbb{Z}^d} |f(\tfrac{z}{N})| \right) \left( \frac{\sqrt{h} \vee |\tfrac{x}{N} - \tfrac{y}{N}|}{\sqrt{t-h}} \right)^\gamma .$$

Hence, by taking limits as in (c.8), the first term in (c.9) vanishes.

For the second term in (c.9), we apply the kernel upper bound (c.1) to obtain

$$\left| S^N_{hN^2,tN^2} f(\tfrac{y}{N}) - S^N_{0,hN^2} S^N_{hN^2,tN^2} f(\tfrac{y}{N}) \right|$$
$$\leq \frac{C}{N^d(\tfrac{1}{N} \vee \sqrt{h})^d} \sum_{\tfrac{z}{N} \in \tfrac{\mathbb{Z}^d}{N}} e^{-\frac{|\tfrac{y}{N} - \tfrac{z}{N}|}{\tfrac{1}{N} \vee \sqrt{h}}} \left| S^N_{hN^2,tN^2} f(\tfrac{y}{N}) - S^N_{hN^2,tN^2} f(\tfrac{z}{N}) \right| . \quad \text{(c.10)}$$

By applying, then, Nash continuity estimate (c.2), we further get

$$\left| S^N_{hN^2,tN^2} f(\tfrac{y}{N}) - S^N_{0,hN^2} S^N_{hN^2,tN^2} f(\tfrac{y}{N}) \right|$$
$$\leq \frac{c \sup_{u \in \mathbb{R}^d} |f(u)| C}{(t-h)^{\gamma/2}} \left\{ \frac{1}{N^d(\tfrac{1}{N} \vee \sqrt{h})^d} \sum_{\tfrac{z}{N} \in \tfrac{\mathbb{Z}^d}{N}} e^{-\frac{|\tfrac{y}{N} - \tfrac{z}{N}|}{\tfrac{1}{N} \vee \sqrt{h}}} \left| \tfrac{y}{N} - \tfrac{z}{N} \right|^\gamma \right\} . \quad \text{(c.11)}$$

We observe that, because we will take the limit $h \to 0$, we can bound above the term outside the curly brackets in (c.11) by a positive constant independent of $N$ and $h$. Concerning the expression between curly brackets in (c.11), by the change of variables $\frac{x}{N(\tfrac{1}{N} \vee \sqrt{h})} \mapsto w$, we obtain

$$\frac{1}{N^d(\tfrac{1}{N} \vee \sqrt{h})^d} \sum_{\tfrac{z}{N} \in \tfrac{\mathbb{Z}^d}{N}} e^{-\frac{|\tfrac{y}{N} - \tfrac{z}{N}|}{\tfrac{1}{N} \vee \sqrt{h}}} \left| \tfrac{y}{N} - \tfrac{z}{N} \right|^\gamma \leq (\tfrac{1}{N} \vee \sqrt{h})^\gamma \, \mathcal{I}_N ,$$

where $\mathcal{I}_N$ is non-negative, is given by

$$\mathcal{I}_N := \sup_{M \in [1, N\sqrt{\theta}]} \frac{1}{M^d} \sum_{w \in \frac{\mathbb{Z}^d}{M}} e^{-|w|} |w|^\gamma$$

and is such that $\limsup_{N \to \infty} \mathcal{I}_N < \infty$. Therefore, we have

$$\lim_{\theta \to 0} \limsup_{N \to \infty} \sup_{h \in [0,\theta]} \left| S^N_{hN^2, tN^2} f(\tfrac{y}{N}) - S^N_{0, hN^2} S^N_{hN^2, tN^2} f(\tfrac{y}{N}) \right| \leq \lim_{\theta \to 0} \theta^{\gamma/2} \limsup_{N \to \infty} \mathcal{I}_N = 0.$$

This concludes the proof of (c.8) for the case $n = 1$.

We then consider the case $n = 2$:

$$\left| \mathbb{E}\left[ f_1\left( \frac{X^x_{hN^2, (t_1+h)N^2}}{N} \right) f_2\left( \frac{X^x_{hN^2, (t_2+h)N^2}}{N} \right) \right] - \mathbb{E}\left[ f_1\left( \frac{X^y_{0, t_1 N^2}}{N} \right) f_2\left( \frac{X^y_{0, t_2 N^2}}{N} \right) \right] \right|,$$

which, rewritten in terms of the random walk semigroups, reads

$$\left| S^N_{hN^2, (t_1+h)N^2} \left( f_1 \cdot S^N_{(t_1+h)N^2, (t_2+h)N^2} f_2 \right)(\tfrac{x}{N}) - S^N_{0, t_1 N^2} \left( f_1 \cdot S^N_{t_1 N^2, t_2 N^2} f_2 \right)(\tfrac{y}{N}) \right|. \quad (\text{c.12})$$

We bound from above (c.12) by

$$\left| S^N_{hN^2, (t_1+h)N^2} \left( f_1 \cdot S^N_{(t_1+h)N^2, (t_2+h)N^2} f_2 \right)(\tfrac{x}{N}) - S^N_{0, t_1 N^2} \left( f_1 \cdot S^N_{(t_1+h)N^2, (t_2+h)N^2} f_2 \right)(\tfrac{y}{N}) \right|$$
$$+ \left| S^N_{0, t_1 N^2} \left( f_1 \cdot S^N_{(t_1+h)N^2, (t_2+h)N^2} f_2 \right)(\tfrac{y}{N}) - S^N_{0, t_1 N^2} \left( f_1 \cdot S^N_{t_1 N^2, t_2 N^2} f_2 \right)(\tfrac{y}{N}) \right|. \quad (\text{c.13})$$

We first observe that, because

$$\sup_{x \in \mathbb{Z}^d} |f_1(\tfrac{x}{N}) S^N_{(t_1+h)N^2, (t_2+h)N^2} f_2(\tfrac{x}{N})| \leq \left( \sup_{u \in \mathbb{R}^d} |f_1(u)| \right) \left( \sup_{u \in \mathbb{R}^d} |f_2(u)| \right),$$

i.e. the supremum of the function $f_1 \cdot S^N_{hN^2, (t_1+h)N^2} f_2$ is independent of $N \in \mathbb{N}$, $t_1, t_2$ and $h \geq 0$, the same argument as for the case $n = 1$ ensures that, by taking limits as in (c.8), the first term in (c.13) vanishes. Secondly, we have that

$$\left| S^N_{0, t_1 N^2} \left( f_1 \cdot S^N_{(t_1+h)N^2, (t_2+h)N^2} f_2 \right)(\tfrac{y}{N}) - S^N_{0, t_1 N^2} \left( f_1 \cdot S^N_{t_1 N^2, t_2 N^2} f_2 \right)(\tfrac{y}{N}) \right|$$
$$\leq \sup_{z \in \mathbb{Z}^d} |f_1(\tfrac{z}{N})| \sup_{z \in \mathbb{Z}^d} \left| S^N_{(t_1+h)N^2, (t_2+h)N^2} f_2(\tfrac{z}{N}) - S^N_{t_1 N^2, t_2 N^2} f_2(\tfrac{z}{N}) \right|. \quad (\text{c.14})$$

Therefore, by taking limits as in (c.8), we can show – as in the first part of the proof, case $n = 1$ – that the r.h.s. in (c.14) and, in turn, the second term in (c.13) vanish.

By repeatedly applying the Markov property, similar arguments yield (c.8) for all $n \in \mathbb{N}$. □

**Tightness.** We end this section with the proof of tightness of the sequence of rescaled random walks

$$\left\{ \frac{X^{x_N}_{0,tN^2}}{N} : t \in [0, T] \right\} \tag{C.15}$$

with arbitrary starting points. Because tightness of these random walks as elements in $\mathsf{D}([0, T], \mathbb{R}^d)$ is implied by tightness of the processes

$$\left\{ \left| \frac{X^{x_N}_{0,tN^2}}{N} \right| : t \in [0, T] \right\} \tag{C.16}$$

in $\mathsf{D}([0, T], \mathbb{R})$ (with "$|\cdot|$" denoting Euclidean distance in $\mathbb{R}^d$), we adopt the tightness criterion presented in Theorem B.4 in Appendix B. The application of this criterion allows us to avoid estimates on certain hitting times as those appearing, for instance, in [2, Proposition 5.13] or [10, Proposition 3.10]. Moreover, the proof that we present relies solely on the kernel upper bound (C.1), item (I) in Proposition C.1, and not on the second bound – item (II) – in Proposition C.1.

**Lemma C.3** (TIGHTNESS). *For all $T > 0$, $u \in \mathbb{R}^d$ and $\{x_N : N \in \mathbb{N}\} \subset \mathbb{Z}^d$ such that $\frac{x_N}{N} \to u$ as $N \to \infty$, the family of processes given as in (C.15) is tight in $\mathsf{D}([0, T], \mathbb{R}^d)$.*

*Proof.* For condition (T1*), we note that there exists an $N_* \in \mathbb{N}$ such that, for all $N \geq N_*$, $\left| \frac{x_N}{N} - u \right| < 1$. Then, for all $m, \ell > 0$ with $m > (|u| + 1) + \ell$, we have

$$\mathbb{P}\left( \left| \frac{X^{x_N}_{0,tN^2}}{N} \right| > m \right) \leq \mathbb{P}\left( \left| \frac{X^{x_N}_{0,tN^2}}{N} - \frac{x_N}{N} \right| > \ell \right),$$

and, thus,

$$0 \leq \lim_{m \to \infty} \limsup_{N \to \infty} \mathbb{P}\left( \left| \frac{X^{x_N}_{0,tN^2}}{N} \right| > m \right) \leq \lim_{\ell \to \infty} \limsup_{N \to \infty} \mathbb{P}\left( \left| \frac{X^{x_N}_{0,tN^2}}{N} - \frac{x_N}{N} \right| > \ell \right).$$

Therefore, by the kernel upper bound (C.1), we obtain:

$$\lim_{\ell \to \infty} \limsup_{N \to \infty} \mathbb{P}\left( \left| \frac{X^{x_N}_{0,tN^2}}{N} - \frac{x_N}{N} \right| > \ell \right) \leq \lim_{\ell \to \infty} \limsup_{N \to \infty} \frac{C}{N^d (\frac{1}{N} \vee \sqrt{t})^d} \sum_{|\frac{y}{N}| > \ell} e^{-\frac{|\frac{y}{N}|}{\frac{1}{N} \vee \sqrt{t}}}$$

$$= \lim_{\ell \to \infty} \frac{C}{t^{\frac{d}{2}}} \int_{|v| > \ell} e^{-\frac{|v|}{\sqrt{t}}} \mathrm{d}v = 0.$$

For condition (T2*), we verify, for all $\varepsilon > 0$, conditions (i)–(v) in Theorem B.4. Let us fix $\varepsilon > 0$ and choose $h_\varepsilon > 0$ and $N_\varepsilon \in \mathbb{N}$ such that

$$\frac{1}{N_\varepsilon} \leq \sqrt{h_\varepsilon} \leq \frac{\varepsilon}{d + \frac{\sqrt{d}}{2} + 1} \,. \tag{C.17}$$

We define, for all $N \geq N_\varepsilon$ and $h \in (0, h_\varepsilon]$,

$$\psi_\varepsilon^N(h) := \frac{C}{N^d (\frac{1}{N} \vee \sqrt{h})^d} \sum_{|\frac{y}{N}| > \varepsilon} e^{-\frac{|\frac{y}{N}|}{\frac{1}{N} \vee \sqrt{h}}} \,, \tag{C.18}$$

with $\psi_\varepsilon^N(0) := 0$. Due to the kernel upper bound (C.1), we get (ii), i.e., for all $h \in [0, h_\varepsilon]$ and $t \in [0, T - h]$,

$$\mathbb{P}\left( \left| \frac{X_{0,(t+h)N^2}^{x_N}}{N} - \frac{X_{0,tN^2}^{x_N}}{N} \right| > \varepsilon \,\middle|\, \left| \frac{X_{0,tN^2}^{x_N}}{N} \right| \right) \leq \psi_\varepsilon^N(h), \quad \mathbb{P}\text{-}a.s. \,.$$

Moreover, because $\sqrt{h_\varepsilon} \leq \frac{\varepsilon}{d + \frac{\sqrt{d}}{2} + 1} < \frac{\varepsilon}{d}$ in (C.17), the function $\psi_\varepsilon^N : [0, h_\varepsilon] \to [0, \infty)$ in (C.18) is non-decreasing. Indeed, for all $h > 0$ such that $\sqrt{h} < \frac{\varepsilon}{d}$,

$$\frac{d}{dh}\left( \frac{1}{h^{d/2}} \sum_{|\frac{x}{N}| > \varepsilon} e^{-\frac{|\frac{x}{N}|}{\sqrt{h}}} \right) = \frac{1}{2 h^{(d+3)/2}} \sum_{|\frac{x}{N}| > \varepsilon} e^{-\frac{|\frac{x}{N}|}{\sqrt{h}}} \left( |\tfrac{x}{N}| - d \sqrt{h} \right) > 0 \,. \tag{C.19}$$

We further introduce the following functions:

$$\psi_\varepsilon(h) := \frac{C}{h^{d/2}} \int_{|v| > \varepsilon} e^{-\frac{|v|}{\sqrt{h}}} dv \tag{C.20}$$

for $h \in (0, h_\varepsilon]$ and $\psi_\varepsilon(0) := 0$ for which $\psi_\varepsilon(h) \to 0$ as $h \to 0$, i.e. (v), holds, as well as

$$\phi_\varepsilon^N := \phi_{\varepsilon,1}^N + \phi_{\varepsilon,2}^N := \sup_{h \in (0, \frac{1}{N^2})} \left| \psi_\varepsilon(h) - \psi_\varepsilon(\tfrac{1}{N^2}) \right| + \sup_{h \in [\frac{1}{N^2}, h_\varepsilon]} \left| \psi_\varepsilon(h) - \psi_\varepsilon^N(h) \right| \,;$$

where

$$\phi_{\varepsilon,1}^N := \sup_{h \in (0, \frac{1}{N^2})} \left| \frac{C}{h^{d/2}} \int_{|v| > \varepsilon} e^{-\frac{|v|}{\sqrt{h}}} dv - \frac{C}{\frac{1}{N^d}} \int_{|v| > \varepsilon} e^{-\frac{|v|}{\frac{1}{N}}} dv \right|$$

$$\phi_{\varepsilon,2}^N := \sup_{h \in [\frac{1}{N^2}, h_\varepsilon]} \left| \frac{C}{h^{d/2}} \int_{|v| > \varepsilon} e^{-\frac{|v|}{\sqrt{h}}} dv - \frac{C}{N^d h^{d/2}} \sum_{|\frac{y}{N}| > \varepsilon} e^{-\frac{|\frac{y}{N}|}{\sqrt{h}}} \right| \,.$$

As a consequence of these definitions, for all $N \geq N_\varepsilon$ and $h \in [0, h_\varepsilon]$, we have (iii), i.e.

$$\psi_\varepsilon^N(h) \leq \psi_\varepsilon(h) + \phi_\varepsilon^N . \qquad \text{(C.21)}$$

Indeed, for all $N \geq N_\varepsilon$ and $h \in [0, h_\varepsilon]$,

$$\begin{aligned}
\psi_\varepsilon^N(h) &\leq \psi_\varepsilon(h) + |\psi_\varepsilon(h) - \psi_\varepsilon^N(h)| \\
&\leq \psi_\varepsilon(h) + |\psi_\varepsilon(h) - \psi_\varepsilon(h \vee \tfrac{1}{N^2})| + |\psi_\varepsilon(h \vee \tfrac{1}{N^2}) - \psi_\varepsilon^N(h)| \\
&\leq \psi_\varepsilon(h) + \sup_{h \in (0, \frac{1}{N^2})} |\psi_\varepsilon(h) - \psi_\varepsilon(\tfrac{1}{N^2})| + \sup_{h \in [\frac{1}{N^2}, h_\varepsilon]} |\psi_\varepsilon(h) - \psi_\varepsilon^N(h)| .
\end{aligned}$$

We claim that $\phi_\varepsilon^N \to 0$ as $N \to \infty$, i.e. (iv). Indeed, by the definition of $\phi_{\varepsilon,1}^N$ and $\psi_\varepsilon(h) \to 0$ as $h \to 0$, we have

$$0 \leq \phi_{\varepsilon,1}^N \xrightarrow[N \to \infty]{} 0 . \qquad \text{(C.22)}$$

Furthermore, by viewing the (improper) Riemann integral $\psi_\varepsilon(h)$ as the limit of Riemann sums over partitions consisting of hypercubes of size $\frac{1}{N}$ centered at $\frac{y}{N} \in \frac{\mathbb{Z}^d}{N} \subset \mathbb{R}^d$ and by monotonicity of Riemann upper and lower sums, we have

$$\begin{aligned}
0 \leq \phi_{\varepsilon,2}^N &\leq \sup_{h \in [\frac{1}{N^2}, h_\varepsilon]} \frac{C}{N^d h^{d/2}} \sum_{|\frac{y}{N}| > \varepsilon} \left( e^{-\frac{|\frac{y}{N}| - \frac{\sqrt{d}}{2N}}{\sqrt{h}}} - e^{-\frac{|\frac{y}{N}| + \frac{\sqrt{d}}{2N}}{\sqrt{h}}} \right) \\
&\leq \sup_{h \in [\frac{1}{N^2}, h_\varepsilon]} \frac{\sqrt{d}}{N \sqrt{h}} \left( \frac{C}{N^d h^{d/2}} \sum_{|\frac{y}{N}| > \varepsilon} e^{-\frac{|\frac{y}{N}| - \frac{\sqrt{d}}{2N}}{\sqrt{h}}} \right) =: \sup_{h \in [\frac{1}{N^2}, h_\varepsilon]} \widehat{\psi}_\varepsilon^N(h) .
\end{aligned}$$

Positivity of the derivative of $\widehat{\psi}_\varepsilon^N$ w.r.t. $h \in [0, h_\varepsilon]$ (cf. (C.19) for an analogous computation) due to the upper bound (C.17) on $\sqrt{h_\varepsilon}$ ensures that $\widehat{\psi}_\varepsilon^N$ is non-decreasing as a function of $h \in [\frac{1}{N^2}, h_\varepsilon]$. Therefore,

$$0 \leq \phi_{\varepsilon,2}^N \leq \widehat{\psi}_\varepsilon^N(h_\varepsilon) \leq \frac{\sqrt{d}}{N \sqrt{h_\varepsilon}} e^{\frac{\sqrt{d}}{2}} \psi_\varepsilon^N(h_\varepsilon) \xrightarrow[N \to \infty]{} 0 . \qquad \text{(C.23)}$$

□

**Acknowledgments.** We warmly thank S.R.S. Varadhan for many enlightening discussions at an early stage of this work. We are indebted to Francesca Collet for fruitful discussions and constant support all throughout this work. We thank Simone Floreani and Alberto Chiarini for helpful conversations on the final part of this paper as well as both referees for their careful reading and for raising relevant issues on some weak points contained in a previous version of this manuscript; we believe this helped us to improve it.

Part of this work was done during the authors' stay at the Institut Henri Poincaré

(UMS 5208 CNRS-Sorbonne Université) – Centre Emile Borel during the trimester *Stochastic Dynamics Out of Equilibrium*. The authors thank this institution for hospitality and support (through LabEx CARMIN, ANR-10-LABX-59- 01). F.S. thanks laboratoire MAP5 of Université de Paris, and E.S. thanks Delft University, for financial support and hospitality. F.S. acknowledges NWO for financial support via the TOP1 grant 613.001.552 as well as funding from the European Union's Horizon 2020 research and innovation programme under the Marie-Skłodowska-Curie grant agreement No. 754411. This research has been conducted within the FP2M federation (CNRS FR 2036).